\documentclass[envcountsame]{llncs}
\usepackage{latexsym}

\title{\bf The Lov\'asz Number of Random Graphs}

\author{Amin Coja-Oghlan\thanks{Research supported by the Deutsche Forschungsgemeinschaft
	(grant DFG FOR 413/1-1).\hspace{5cm} Some of the material has appeared in a preliminary version
	in Proc.\ STACS 2003 and Proc.\ RANDOM 2003}}
\date{\today}
\institute{Humboldt-Universit\"{a}t zu Berlin, Institut f\"{u}r Informatik,\\Unter den Linden 6, 10099 Berlin, Germany
\\\email{coja@informatik.hu-berlin.de}}

\newcommand\frakfamily{\usefont{U}{euf}{m}{n}} %eufm10
\DeclareTextFontCommand{\frak}{\frakfamily}
\newcommand\msymbfamily{\usefont{U}{msb}{m}{n}}
\DeclareTextFontCommand{\msymb}{\msymbfamily}

\newcommand\myhomepage{http://www.informatik.hu-berlin.de/$\sim$coja/}

\newcommand\eul{\mathrm{e}}
\newcommand\eps{\varepsilon}
\newcommand\Alg{\mathcal{A}}

\newcommand\Erw{\mathrm{E}}

\newcommand\pr{\mathrm{P}}
\newcommand\gnr{G_{n,r}}
\newcommand\gnp{G_{n,p}}

\newcommand\conf{\mathcal{C}}

\newcommand{\bink}[2]
    {{{#1}\choose {#2}}}

\newcommand\NP{\mathrm{NP}}

\newcommand\thet{\vartheta}
\newcommand\thetI{\vartheta_{1/2}}
\newcommand\thetII{\vartheta_2}
\newcommand\bthet{\bar\thet}
\newcommand\bthetI{\bar\thetI}
\newcommand\bthetII{\bar\thetII}
\newcommand\vchi{\bthetI}

\newcommand\bc[1]{\left({#1}\right)}
\newcommand\bcfr[2]{\bc{\frac{#1}{#2}}}
\newcommand\RR{\mathbf{R}}

\newcommand{\bck}[1]{\langle #1\rangle}
\newcommand\scal[2]{\bck{{#1},{#2}}}

\newcommand{\vecone}{\vec{1}}
\newcommand\SMC{\mathrm{SMC}}

\newcommand\SDP{\mathrm{SDP}}
\newcommand\MKC{\textrm{MAX $k$-CUT}}
\newcommand\mkc{\mathrm{MC}_k}

\newcommand{\Lovasz}{Lov\'asz}
\newcommand{\Juhasz}{Juh\'asz}

\begin{document}

\maketitle

\spnewtheorem{Algo}[theorem]{Algorithm}{\bfseries}{}

\begin{abstract}
We study the \Lovasz\ number $\thet$
along with two further SDP relaxations $\thetI$, $\thetII$
of the independence number and the corresponding relaxations
$\bthet$, $\bthetI$, $\bthetII$ of the chromatic number on random graphs $\gnp$.
We prove that $\thet,\thetI,\thetII(\gnp)$ are concentrated
about their means, and that $\bthet,\bthetI,\bthetII(\gnp)$
in the case $p<n^{-1/2-\eps}$ are concentrated in intervals of constant length.
Moreover, extending a result of \Juhasz\ \cite{Juhasz}, we show that
$\thet,\thetI,\thetII(\gnp)=\Theta(\sqrt{n/p})$
and that $\bthet,\bthetI,\bthetII(\gnp)=\Theta(\sqrt{np})$ for $c_0/n\leq p\leq1/2$.
As an application, we give an improved algorithm for
approximating the independence number of $\gnp$ in polynomial expected time, thereby
extending a result of Krivelevich and Vu \cite{KrivVu}.
We also improve on the analysis of an algorithm of Krivelevich \cite{KrivDecide}
for deciding whether $\gnp$ is $k$-colorable.\\
\emph{Topics and key words: \Lovasz\ number, vector chromatic number,
	random graphs, maximum independent set problem, graph coloring}
\end{abstract}

\section{Introduction and Results}

Given a graph $G=(V,E)$,
let $\alpha(G)$ be the independence number, let
$\omega(G)$ be the clique number, and let $\chi(G)$ be the chromatic number of $G$.
Further, let $\bar G$ signify the complement of $G$.
Since it is NP-hard to compute any of $\alpha(G)$, $\omega(G)$ or $\chi(G)$,
it is remarkable that there exists an efficiently computable function $\thet(G)$
that is ``sandwiched'' between $\alpha(G)$ and $\chi(\bar G)$, i.e.\
	$\alpha(G)\leq\thet(G)\leq\chi(\bar G).$
Passing to complements, and letting $\bthet(G)=\thet(\bar G)$, we have
	$\omega(G)\leq\bthet(G)\leq\chi(G)$.
The function $\thet$ was introduced by \Lovasz\  \cite{Lovasz},
and is called the \emph{\Lovasz\  number} of $G$ (cf.\ also \cite{Knuth}).
The \Lovasz\  number can be seen as a semidefinite programming (``SDP'')
relaxation of the independence number,
and is therefore comptable in polynomial time within any precision \cite{GLS}.

Though $\thet(G)$ is sandwiched between $\alpha(G)$ and $\chi(\bar G)$,
Feige \cite{FeigeTheta} proved that
the gap between $\alpha(G)$ and $\thet(G)$ or between $\chi(\bar G)$ and
$\thet(G)$ can be as large as $n^{1-\eps}$, $\eps>0$.
Indeed, unless NP$=$coRP, none of $\alpha(G)$, $\omega(G)$,
$\chi(G)$ can be approximated within a factor of $n^{1-\eps}$, $\eps>0$, in polynomial time
\cite{Hastad,Zeroknowledge}.
However, though there exist graphs $G$ such that $\thet(G)$ is not a good
approximation of $\alpha(G)$ (or $\bthet(G)$ of $\chi(G)$), it might be the case
that the \Lovasz\  number performs well on ``average'' instances.
In fact, several algorithms for random and semirandom graph problems are
based on computing $\thet$ \cite{Coja,CojaTaraz,FindLarge,Feige,FeigeKraut}.
Therefore, the aim of this paper is to study the \Lovasz\  number of random graphs
more thoroughly.

The standard model of a random graph is the binomial model $\gnp$, pioneered by Erd\H{o}s and Renyi.
We let $0<p=p(n)<1$ be a number that may depend on $n$.
Let $V=\{1,\ldots,n\}$.
Then the random graph $\gnp$ is obtained by
including each of the $\bink{n}{2}$ possible edges $\{v,w\}$, $v,w\in V$,
with probability $p$ independently.
Though $\gnp$ may fail to model some types of input instances appropriately,
both the combinatorial structure and the algorithmic theory of $\gnp$
are of fundamental interest \cite{BB,JLR,FMcD}.
We say that $\gnp$ has some property $A$ \emph{with high probability} (whp.), if
$\lim_{n\rightarrow\infty}\pr(\gnp\textrm{ has property }A)=1$.

In addition to the \Lovasz\ number, we also address two further natural
SDP relaxations $\thetI$, $\thetII$ of $\alpha$ (cf.\ \cite{Szegedy}) on random graphs.
These relaxations satisfy
	$\alpha(G)\leq\thetI(G)\leq\thet(G)\leq\thetII(G)\leq\chi(\bar G),$
for all $G$, i.e.\ $\thetI$ is the strongest relaxation of $\alpha$ among $\thetI,\thet,\thetII$.
Passing to complements, and setting $\bthet_i(G)=\thet_i(\bar G)$ ($i=1/2,2$),
one gets $\omega(G)\leq\bthetI(G)\leq\bthet(G)\leq\bthetII(G)\leq\chi(G)$,
i.e.\ $\bthetII$ is the strongest relaxation of $\chi$.
The relaxation $\bthetI(G)$ coincides with the well-known \emph{vector chromatic number}
$\vec\chi(G)$ of Karger, Motwani, and Sudan \cite{KMS}.

\subsubsection{The concentration of $\thet$, $\bthet$, etc.}
Facing a real-valued random variable $X(\gnp)$, there are two obvious questions
to ask.
\begin{enumerate}
\item What is the mean of $X(\gnp)$?
\item Is $X(\gnp)$ concentrated about its mean?
\end{enumerate}
The main contributions of this paper are concentration results on $\thet$, $\bthet$ etc.,
i.e.\ concern the second question.
Such results are important for instance in the design of algorithms with
a polynomial expected running time.
First, we show that the probability that $\thet(\gnp)$, $\thetI(\gnp)$, or $\thetII(\gnp)$ is
far from its median is exponentially small.

\begin{theorem}\label{Thm_Largedev}
Suppose that $p\leq 0.99$, and that $n\geq n_0$ for a certain constant $n_0>0$.
Let $m$ be a median of $\thet(\gnp)$, i.e.
$\pr(\thet(\gnp)\leq m)\geq 1/2$ and $\pr(\thet(\gnp)\geq m)\geq 1/2$.
\begin{enumerate}
\item[i.]  Let $\xi\geq\max\{10,m^{1/2}\}$.
	Then
	 $\pr(\thet(\gnp)\geq m+\xi)\leq 30\exp(-\xi^2/(5m+10\xi)). $
\item[ii.] Let $\xi>10$. Then
	$\pr(\thet(\gnp)\leq m-\xi)\leq3\exp(-\xi^2/10m).$
\end{enumerate}
The same holds with $\thet$ replaced by $\thetI$ or by $\thetII$.
\end{theorem}

Up to the constants involved, the right hand sides of the bounds in Thm.\ \ref{Thm_Largedev}
are similar to well-known bounds on the tails of the binomial distribution (e.g.\ \cite[p.\ 26]{JLR}).
The proof of Thm.\ \ref{Thm_Largedev} is based on Talagrand's inequality.
%, see Sec.\ \ref{Sec_Largedev} for details.

A remarkable fact concerning the chromatic number of sparse random graphs $\gnp$,
$p\leq n^{-\eps-1/2}$, is that $\chi(\gnp)$ is concentrated in an interval of constant length.
Indeed, Shamir and Spencer \cite{ShamirSpencer} proved that there is a function $u=u(n,p)$
such that in the case $p=n^{-\beta}$, $1/2<\beta<1$, we have
	$\pr(u\leq\chi(\gnp)\leq u+\lceil(2\beta+1)/(2\beta-1)\rceil)=1-o(1).$
Furthermore, \L uczak \cite{Luczak} showed that in the case $5/6<\beta<1$, the chromatic number
is concentrated in width one, which is best possible.
In fact, Alon and Krivelevich \cite{AlonKriv} could prove that two point concentration
holds for the entire range $p=n^{-\beta}$, $1/2<\beta<1$.
The two following theorems state similar results as given by Shamir and Spencer
and by \L uczak for the relaxations $\vchi(\gnp)$, $\bthet(\gnp)$, and $\bthetII(\gnp)$
of the chromatic number.

\begin{theorem}\label{Thm_ShamirSpencer}
Suppose that $c_0/n\leq p\leq n^{-\beta}$ for some large constant $c_0>0$ and
some number $1/2<\beta<1$.
Then $\vchi(\gnp)$, $\bthet(\gnp)$, $\bthetII(\gnp)$ are concentrated in width
$s=\frac{2}{2\beta-1}+o(1)$,
i.e.\ there exist numbers $u$, $u'$, $u''$ depending on $n$ and $p$ such that whp.\
	$$u\leq\vchi(\gnp)\leq u+s,\ u'\leq\bthet(\gnp)\leq u'+s,\textrm{ and }u''\leq\vchi(\gnp)\leq u''+s.$$
\end{theorem}

\begin{theorem}\label{Thm_Luczak}
Suppose that $c_0/n<p\leq n^{-5/6-\delta}$ for some large constant $c_0$
and some $\delta>0$.
Then $\vchi(\gnp)$, $\bthet(\gnp)$, and $\bthetII(\gnp)$ are concentrated in width $1$.
\end{theorem}

In contrast to the chromatic number, $\vchi$, $\bthet$, and $\bthetII$ need not be integral.
Therefore, the above results do \emph{not} imply that
$\vchi(\gnp)$, $\bthet(\gnp)$, $\bthetII(\gnp)$ are concentrated on a constant number of points.
The proofs of Thms.\ \ref{Thm_ShamirSpencer} and \ref{Thm_Luczak}
are given in Sec.\ \ref{Sec_Sharp}.

\subsubsection{The probable value of $\thet(\gnp)$, $\bthet(\gnp)$, etc.}
Concerning the probable value of $\thet(\gnp)$ and $\bthet(\gnp)$,
\Juhasz\ \cite{Juhasz} gave the following partial answer:
If $\ln(n)^6/n\ll p\leq 1/2$, then with high probability we have
	$\thet(\gnp)=\Theta(\sqrt{n/p})\textrm{ and }\bthet(\gnp)=\Theta(\sqrt{np}).$
However, we shall indicate in Sec.~\ref{Sec_thetgnp} that \Juhasz's proof
fails in the case of sparse random graphs (e.g.\ $np=O(1)$).
Making use of the above concentration results on $\thet$, $\bthet$ etc.,
we can compute the probable value not only of $\thet(\gnp)$ and $\bthet(\gnp)$, but
also of $\thet_i(\gnp)$ and $\bthet_i(\gnp)$, $i=1/2,2$, for essentially the
entire range of edge probabilities $p$.
To the best of the author's knowledge, no previous results concerning
$\thet_i(\gnp)$ and $\bthet_i(\gnp)$, $i=1/2,2$, occur in the literature.
Note that we only need to consider edge probabilities $p\leq1/2$, because
$G_{n,1-p}=\bar\gnp$.

\begin{theorem}\label{Thm_thetgnp}
Suppose that $c_0/n\leq p\leq1/2$ for some large constant $c_0>0$.
Then there exist constants $c_1,c_2,c_3,c_4>0$ such that
	\begin{eqnarray}
	&c_1\sqrt{n/p}\leq\thetI(\gnp)\leq\thet(\gnp)\leq\thetII(\gnp)\leq c_2\sqrt{n/p}&\label{eqthetgnp}\\
	\textrm{and}&c_3\sqrt{np}\leq\bthetI(\gnp)\leq\bthet(\gnp)\leq\bthetII(\gnp)\leq c_4\sqrt{np}&\nonumber
	\end{eqnarray}
with high probability. More precisely,
	\begin{equation}\label{eqexp}
	\pr(c_3\sqrt{np}\leq\vchi(\gnp)\leq\bthet(\gnp)\leq\bthetII(\gnp))\geq1-\exp(-n).
	\end{equation}
\end{theorem}

Assume that $c_0/n\leq p=o(1)$.
Then $\alpha(\gnp)\sim2\ln(np)/p$ and $\chi(\gnp)\sim np/(2\ln(np))$ whp.\  (cf.\ \cite{JLR}).
Hence, Thm.\ \ref{Thm_thetgnp} shows that $\thetII(\gnp)$ (resp.\ $\bthetI(\gnp)$)
approximates $\alpha(\gnp)$ (resp.\ $\chi(\gnp)$) within a factor of $O(\sqrt{np})$.
In fact, if $np=O(1)$, then we get a constant factor approximation.
On the other hand, as $\alpha(G_{n,1/2})\sim 2\log_2(n)$ and $\chi(G_{n,1/2})\sim n/(2\log_2(n))$,
in the random graph $G=G_{n,1/2}$
the gap between $\thetI(G)$ (resp.\ $\bthetII(G)$) and $\alpha(G)$ (resp.\ $\chi(G)$)
is as large as $n^{1/2-\eps}$ whp.
Our estimate on the probable value of the vector chromatic number $\bthetI(\gnp)$
in Thm.~\ref{Thm_thetgnp} answers a question of Krivelevich~\cite{KrivDecide}.

As a consequence of the upper bound on $\bthetII(\gnp)$ in Thm.\ \ref{Thm_thetgnp},
we obtain a lower bound on the probable value of the SDP relaxation $\SDP_k$ of \MKC\
due to Frieze and Jerrum \cite{FJ}.

\begin{corollary}\label{Cor_SDPk}
Let $k\geq2$ be an integer.
Suppose that $c_0k^2/n\leq p\leq 1/2$.
Then whp.\ we have
	$\SDP_k(\gnp)\geq\bc{1-\frac{1}{k}}\bink{n}{2}p+c_1n^{3/2}p^{1/2}$
for some constant $c_1>0$.
\end{corollary}

Corollary \ref{Cor_SDPk} complements an upper bound on $\SDP_k(\gnp)$ due to Coja-Oghlan,
Moore, and Sanwalani \cite{MAXCUT}, who proved that there is a constant $c_2>0$
such that
	$\SDP_k(\gnp)\leq(1-1/k)\bink{n}{2}p+c_2n^{3/2}p^{1/2}$
whp.
In contrast, the weight $\mkc(\gnp)$ of a \MKC\ of $\gnp$ is at most
	$$\bc{1-\frac{1}{k}}\bink{n}{2}p+\sqrt{\frac{\ln(k)}{k}}\cdot n^{3/2}p^{1/2}$$
whp.\ (cf.\ \cite{MAXCUT}).
Thus, Cor.\ \ref{Cor_SDPk} shows that for large $k$ there is a moderate gap between
$\SDP_k(\gnp)$ and $\mkc(\gnp)$.

Finally, let us consider the \emph{random regular graph} $\gnr$,
i.e.\ an $r$-regular graph of order $n$ chosen uniformly at random.

\begin{theorem}\label{Thm_ThetReg}
Let $c_0$ be a sufficiently large constant, and let $c_0\leq r=o(n^{1/4})$.
There are constants $c_1,c_2>0$ such that whp.\ the random regular graph $\gnr$
satisfies
	$$c_1n/\sqrt{r}\leq\thetI(\gnr)\leq\thet(\gnr)\leq\thetII(\gnr)\leq c_2n/\sqrt{r}.$$
Moreover, there is a constant $c_3>0$ such that in the case $c_0\leq r=o(n^{1/2})$ we have
 	\begin{equation}\label{eqvchignr}
	\pr(c_3\sqrt{r}\leq\vchi(\gnr)\leq\bthet(\gnr)\leq\bthetII(\gnr))\geq1-\exp(-n).
	\end{equation}
\end{theorem}

\subsubsection{Algorithmic applications.}
There are two types of algorithms for $\NP$-hard random graph problems.
First, there are \emph{heuristics} that \emph{always} run in polynomial time,
and \emph{almost always} output a good solution.
On the other hand, there are algorithms that guarantee some approximation ratio
on \emph{any} input instance, and which have a polynomial \emph{expected}
running time when applied to $\gnp$ (cf.\ \cite{DyerFrieze}).
Here we say that an algorithm $A$ runs in polynomial expected time
if there is a constant $l>0$ such that $\sum_G R_A(G) \,\pr(\gnp=G)=O(n^l)$, where $R_A(G)$
is the running time of $A$ on input $G$
and the sum ranges over all graphs $G$ of order $n$.
In this paper, we are concerned with algorithms with a polynomial expected running time.

First, we consider the maximum independent set problem in random graphs.
Krivelevich and Vu \cite{KrivVu} gave an algorithm
that in the case $p\gg n^{-1/2}$ approximates the independence number of $\gnp$ in polynomial
expected time within a factor of $O(\sqrt{np}/\ln(np))$.
Moreover, they ask whether a similar algorithm exists for smaller values of $p$.
A first answer was obtained by Coja-Oghlan and Taraz \cite{CojaTaraz},
who gave an $O(\sqrt{np}/\ln(np))$-approximative algorithm for the case $p\gg\ln(n)^6/n$.
Using Thms.\ \ref{Thm_Largedev} and \ref{Thm_thetgnp}, we can improve on the analysis given in \cite{CojaTaraz},
thereby answering the question of Krivelevich and Vu in the affirmative.

\begin{theorem}\label{Thm_MIS}
Suppose that $c_0/n\leq p\leq1/2$.
There exists an algorithm $\texttt{ApproxMIS}$ that for any input graph $G$
outputs an independent set of size at least $\alpha(G)\ln(np)/(c_1\sqrt{np})$,
and which applied to $\gnp$ runs in polynomial expected time.
Here $c_0,c_1>0$ denote constants.
\end{theorem}

As a second application, we give an algorithm for deciding within polynomial
expected time whether the input graph is $k$-colorable.
Instead of $\gnp$, we shall even consider the \emph{semirandom model} $\gnp^+$ that allows for an
adversary to add edges to the random graph.
More precisely, the semirandom graph $\gnp^+$ is constructed in two steps as follows.
First, a random graph $G_0=\gnp$ is chosen.
Then, an adversary completes the instance $G=\gnp^+$ by adding arbitrary edges to $G_0$.
We say that \emph{the expected running time of an algorithm $\Alg$ is polynomial over $\gnp^+$},
if there is some constant $l$ such that the expected running time of $\Alg$ is $O(n^l)$
regardless of the behavior of the adversary.

\begin{theorem}\label{Thm_Decide}
Suppose that $k=o(\sqrt{n})$, and that $p\geq c_0k^2/n$, for some constant $c_0>0$.
There exists an algorithm $\texttt{Decide}_k$ that for any input graph $G$ decides whether $G$
is $k$-colorable, and that applied to $\gnp^+$ has a polynomial expected running time.
\end{theorem}

The algorithm $\texttt{Decide}_k$ is essentially identical with
Krivelevich's algorithm for deciding $k$-colorability in polynomial expected time \cite{KrivDecide}.
However, the analysis given in \cite{KrivDecide} requires that $np\geq\exp(\Omega(k))$, whereas
Thm.\ \ref{Thm_Decide} only requires that $np$ is quadratic in $k$.
The improvement results from the fact that the analysis given in this paper relies on
the asymptotics for $\vchi(\gnp)$ derived in Thm.\ \ref{Thm_thetgnp} (instead of
the concept of semi-colorings).
Finally, we prove that our algorithm $\texttt{Decide}_k$ also applies to random regular graphs $\gnr$.

\begin{theorem}\label{Thm_DecideReg}
Suppose that $c_0k^2\leq r=o(n^{1/2})$ for some constant $c_0>0$.
Then, applied to $\gnr$, the algorithm $\texttt{Decide}_k$ has polynomial expected running time.
\end{theorem}

\subsubsection{Organization.}
First we recall the definitions of $\thet$, $\bthet$ etc.\ and prove some elementary
facts in Sec.\ \ref{Sec_Pre}.
Sec.\ \ref{Sec_Conc} deals with the concentration results, and Sec.\ \ref{Sec_thetgnp} contains
the proofs of Thm.\ \ref{Thm_thetgnp} and Cor.\ \ref{Cor_SDPk}.
In Sec.\ \ref{Sec_reg} we prove Thm.\ \ref{Thm_ThetReg},
and finally Sec.\ \ref{Sec_Algo} is devoted to the algorithms.

\subsubsection{Notation.}
Throughout we let $V=\{1,\ldots,n\}$.
If $G=(V,E)$ is a graph, and $U\subset V$, then $N(U)$ is the neighborhood of $U$,
i.e.\ set of all $v\in V$ such that there is $w\in U$ satisfying $\{v,w\}\in E$.
Moreover, $A(G)$ is the adjacency matrix of $G$.
By $\vecone$ we denote the vector with all entries equal to  one in any dimension.
Furthermore, $J$ denotes a square matrix of any size with all entries equal to one.
If $M$ is a real symmetric $n\times n$-matrix, then $\lambda_1(M)\geq\cdots\lambda_n(M)$ signify
the eigenvalues of $M$, and $\|M\|=\max\{\lambda_1(M),-\lambda_n(M)\}$ is the spectral radius of $M$.
We let $\scal{\cdot}{\cdot}$ denote the scalar product of vectors.
By $c_0,c_1,\ldots$ we denote constants, i.e.\ numbers that are independent of $n$ and $p$.

\section{Preliminaries}\label{Sec_Pre}

In this section we recall the definitions of $\thet$, $\thetI$, $\thetII$,
and provide some elementary facts which will be useful later.
We let $G=(V,E)$ be a graph, and let  $\bar G$ be the complement of $G$.
%
%Following \cite{Knuth}, an \emph{orthogonal labeling} of $\bar G$
%is a tuple $(v_1,\ldots,v_n)$ of vectors
%$v_i\in\RR^d$ such that for any two vertices $i,j\in V$, $i\not=j$,
%with $\{i,j\}\in E$ we have  $v_i\perp v_j.$
%Here $d>0$ is any integer.
%Furthermore, the \emph{cost} of a $d$-dimensional vector $a=(a_1,\ldots,a_d)$ is
% $$c(a)=\left\{\begin{array}{cl}a_1^2\|a\|^{-2}&\textrm{ if }a\not=0\\
%        0&\textrm{ otherwise.}\end{array}\right.$$
%Then $0\leq c(a)\leq1$.
%We let
%	$$\thet(G)=\max\left\{\sum_{i=1}^nc(v_i)|\ (v_1,\ldots,v_n)\textrm{ is an orthogonal
%		labeling of }\bar G\right\}\quad\textrm{(cf.\ \cite{Knuth})}.$$
%%%%%%%%%%%%%%%%%%%%%%%%%%%%%%%%%%%
%On the other hand, we can characterize $\thet(G)$ as the solution to an eigenvalue
%minimization problem as follows.
%We call a real symmetric $n\times n$-matrix
%$M=(m_{ij})_{i,j=1,\ldots,n}$
%\emph{$\thet$-feasible} for $G$ if $m_{ii}=1$ for all $i$, and
%$m_{ij}=1$ whenever $\{i,j\}\not\in E$.
%Since $M$ is symmetric, the eigenvalues $\lambda_1(M)\geq\cdots\geq\lambda_n(M)$ are real numbers.
%Then, we have
%	$$\thet(G)=\min\{\lambda_1(M)|\ M\textrm{ is feasible for }G\},$$
%cf.\ \cite{Knuth}.
%%%%%%%%%%%%%%%%%%%%%%%%%%%%%%%%%
%We define $\vchi(G)$ and $\bthetII(G)$ in terms of vector colorings as follows.
Let $(v_1,\ldots,v_n)$ be an $n$-tuple of unit vectors in $\RR^n$, and let $k>1$.
Then $(v_1,\ldots,v_n)$ is a \emph{vector $k$-coloring} of $G$
if $\scal{v_i}{v_j}\leq-1/(k-1)$ for all edges $\{i,j\}\in E$.
Furthermore, $(v_1,\ldots,v_n)$ is a \emph{strict} vector $k$-coloring
if $\scal{v_i}{v_j}=-1/(k-1)$ for all $\{i,j\}\in E$.
Finally, we say that $(v_1,\ldots,v_n)$ is a \emph{rigid} vector $k$-coloring if
$\scal{v_i}{v_j}=-1/(k-1)$ for all $\{i,j\}\in E$ and $\scal{v_i}{v_j}\geq-1/(k-1)$
for all $\{i,j\}\not\in E$.
Following \cite{KMS,GoemansKleinberg,Charikar}, we define
	\begin{eqnarray}
	\vchi(G)&=&\inf\{k>1|\ \textrm{$G$ admits a vector $k$-coloring}\},\nonumber\\
	\bthet(G)=\bthet_1(G)&=&\inf\{k>1|\ \textrm{$G$ admits a strict vector $k$-coloring}\},\label{eqstrictvectorcol}\\
	\bthetII(G)&=&\inf\{k>1|\ \textrm{$G$ admits a rigid vector $k$-coloring}\}.\nonumber
	\end{eqnarray}
Observe that $\vchi(G)$ is precisely the \emph{vector chromatic number}
introduced by Karger, Motwani, and Sudan \cite{KMS};
$\bthetII$ occurs in \cite{GoemansKleinberg,Szegedy}.
Further, we let $\thetI(G)=\bthetI(\bar G)$, $\thet(G)=\thet_1(G)=\bthet(\bar G)$, and
$\thetII(G)=\bthetII(\bar G)$.
%Moreover, letting $\bthet(G)=\thet(\bar G)$, we have
%	\begin{equation}\label{eqstrict}
%	\bthet(G)=\inf\{k>1|\ \textrm{$G$ admits a strict vector $k$-coloring}\},
%	\end{equation}
%cf.\ \cite{KMS}.
It is shown in \cite{KMS} that the above definition of $\thet$ is equivalent
with \Lovasz's original definition \cite{Lovasz}.

\begin{proposition}\label{Prop_facts}
Let $G=(V,E)$ be a graph of order $n$, and let $S\subset V$.
Let $G\lbrack S\rbrack$ denote the subgraph of $G$ induced on $S$.
Then $\thet_i(G)\leq\thet_i(G\lbrack S\rbrack)+\thet_i(G\lbrack V\setminus S\rbrack)$,
	$i\in\{1/2,1,2\}$.
\end{proposition}

Although the result may be known to specialists in the area,
to the best of the author's knowledge it is not explicitly stated (or proved) in the literature.
%Since several proofs in this paper rely on Prop.\ \ref{Prop_facts},
Therefore, we carry out the proof for $\bthetII$; the same argument applies to $\vchi$ and $\bthet$.

\noindent\emph{Proof of Prop.\ \ref{Prop_facts}.}
Let $k>\bthetII(G\lbrack S\rbrack)$ and let $l>\bthetII(G\lbrack V\setminus S\rbrack)$.
Further, let $(a_v)_{v\in S}$ be a rigid vector $k$-coloring of $G\lbrack S\rbrack$,
and let $(b_v)_{v\in V\setminus S}$ be a rigid vector $l$-coloring of
$G\lbrack V\setminus S\rbrack$.
Set
	$$\alpha=\bcfr{l}{(k+l)(k-1)}^{1/2}\textrm{ and }\beta=\bcfr{k}{(k+l)(l-1)}^{1/2}.$$
Embedding the $a_v$'s and $b_w$'s into a high-dimensional space, we may assume that
$a_v\perp b_w$ for all $v\in S$, $w\in V\setminus S$, and that there is a unit vector $z$
such that $z\perp a_v$, $z\perp b_w$ for all $v,w$.
Let
	$$x_v=(1+\alpha^2)^{-1/2}(a_v+\alpha z)\textrm{ and }x_w=(1+\beta^2)^{-1/2}(b_w-\beta z)
		\qquad(v\in S,w\in V\setminus S).$$
Then $\|x_v\|=\|x_w\|=1$.
Moreover, if two vertices $v,v'\in S$ are adjacent, then
	$$\scal{x_v}{x_{v'}}=(1+\alpha^2)^{-1}\bc{\scal{a_v}{a_{v'}}+\alpha^2}=(1+\alpha^2)^{-1}\bc{-\frac{1}{k-1}+\alpha^2}
		=-\frac{1}{k+l-1}.$$
Likewise, if $v,v'\in S$ are non-adjacent, then
	$\scal{x_v}{x_{v'}}\geq-1/(k+l-1)$.
Consequently, $(x_v)_{v\in S}$ is a rigid vector $(k+l)$-coloring of $G\lbrack S\rbrack$.
Similarly, $(x_w)_{w\in V\setminus S}$ is a rigid vector $(k+l)$-coloring of $G\lbrack V\setminus S\rbrack$.
Since $\scal{x_v}{x_w}=-1/(k+l-1)$ for all $v\in S$, $w\in V\setminus S$, $(x_v)_{v\in V}$ is a
rigid vector $(k+l)$-coloring of the entire graph $G$, thereby proving $\bthetII(G)\leq k+l$.
\qed

In addition to Prop.\ \ref{Prop_facts}, we will frequently make use of the well-known
fact that
%	\begin{equation}\label{eqtrivbounds}
$$	\omega(G)\leq\vchi(G)\leq\bthet(G)\leq\bthetII(G)\leq\chi(G)\textrm{ and }
		\alpha(G)\leq\thetI(G)\leq\thet(G)\leq\thetII(G)$$
%	\end{equation}
for all graphs $G=(V,E)$.
The lower bounds $\omega(G)\leq\vchi(G)$, $\alpha(G)\leq\thetI(G)$ are
established in \cite{KMS}.
The upper bound $\bthetII(G)\leq\chi(G)$ can be proved e.g.\ by decomposing $V$
into $\chi(G)$ disjoint independent sets and applying Prop.\ \ref{Prop_facts}.
Moreover, it is obvious from the definitions that for any weak subgraph $H$
of $G$ we have
	\begin{equation}\label{eqmonotone}
	\bthet_i(H)\leq\bthet_i(G)\qquad(i\in\{1/2,1,2\}).
	\end{equation}

In addition to $\thet$, $\thetI$, and $\thetII$, we consider the
following semidefinite relaxation of MAX $k$-CUT, due to Frieze and Jerrum \cite{FJ}.
Let $G$ be a graph with adjacency matrix $A=A(G)=(a_{ij})_{i,j=1,\ldots,n}$, and let $k\geq2$.
Then
    \begin{equation}\label{eqSDPk}
    \SDP_k(G)=\max\sum_{i<j}a_{ij}\frac{k-1}{k}\bc{1-\scal{v_i}{v_j}}\textrm{ s.t.\ }\|v_i\|=1,\ \scal{v_i}{v_j}\geq-\frac{1}{k-1},
    \end{equation}
where the $\max$ is taken over $v_1,\ldots,v_n\in\RR^n$,
is an upper bound on the weight of a MAX $k$-CUT of $G$.
In the case $k=2$, we obtain the semidefinite relaxation $\SMC=\SDP_2$
of MAX CUT invented by Goemans and Williamson \cite{GW}.
In this case, the constraint $\scal{v_i}{v_j}\geq-1/(2-1)=-1$ is void.

\section{The Concentration Results}\label{Sec_Conc}

\subsection{Proof of Theorem \ref{Thm_Largedev}}\label{Sec_Largedev}

\subsubsection{The large deviation result for $\thet$.}

In order to bound the probability that the Lov\'asz number $\thet(\gnp)$
is far from its median, we shall apply the following version of
\emph{Talagrand's inequality} (cf.\ \cite[p. 44]{JLR}).

\begin{theorem}\label{Thm_Talagrand}
Let $\Lambda_1,\ldots,\Lambda_N$ be probability spaces.
Let $\Lambda=\Lambda_1\times\cdots\times\Lambda_N$.
Let $A,B\subset\Lambda$ be measurable sets such that for some $t\geq 0$ the following
condition is satisfied:
For every $b\in B$ there is $\alpha=(\alpha_1,\ldots,\alpha_N)\in\RR^N\setminus\{0\}$
such that for all $a\in A$ we have
 $$\sum_{i:\,a_i\not=b_i}\alpha_i\geq t\bc{\sum_{i=1}^N\alpha_i^2}^{1/2},$$
where $a_i$ (resp.\ $b_i)$ denotes the $i$'th coordinate of $a$ (resp.\ $b$).
Then
 $\pr(A)\pr(B)\leq\exp(-t^2/4).$
\end{theorem}

Let $G=(V,E)$ be a graph.
We need the following equivalent characterization of $\thet(G)$.
A tuple $(v_1,\ldots,v_n)$ of vectors $v_i\in\RR^d$ is called an
\emph{orthogonal labeling} of $\bar G$
if for any two vertices $i,j\in V$, $i\not=j$,
with $\{i,j\}\in E$ we have $v_i\perp v_j$ (cf.\ \cite{Knuth}).
Here $d>0$ is any integer.
Furthermore, the \emph{cost} of a $d$-dimensional vector $a={^t(}a_1,\ldots,a_d)$ is
 $$c(a)=\left\{\begin{array}{cl}a_1^2\|a\|^{-2}&\textrm{ if }a\not=0\\
        0&\textrm{ otherwise.}\end{array}\right.$$
Then $0\leq c(a)\leq1$, and we have
	\begin{equation}\label{eqorth}
	\thet(G)=\max\left\{\sum_{i=1}^nc(v_i)|\ (v_1,\ldots,v_n)\textrm{ is an orthogonal
		labelling of }\bar G\right\}\quad\textrm{(cf.\ \cite{KMS,Knuth})}.
	\end{equation}

The proof of Thm.\ \ref{Thm_Largedev} relies on the following lemma.

\begin{lemma}\label{Lemma_Key}
Let $m$ be a median of $\thet(\gnp)$.
Let $\thet_0>0$ be any number, and let $\xi\geq 10$.
Then
 $$\pr(m+\xi\leq\thet(\gnp)\leq\thet_0)\leq 2\exp(-\xi^2/(5\thet_0)).$$
\end{lemma}
\begin{proof}
For $i\geq 2$, let $\Lambda_i\in\{0,1\}^{i-1}$ consist of the first $i-1$ entries of the $i$th row
of the adjacency matrix of $\gnp$.
Then $\Lambda_2,\ldots,\Lambda_n$ are independent random variables,
and $\Lambda_i$ determines to which of the $i-1$ vertices $1,\ldots,i-1$ vertex $i$
is adjacent.
Therefore, we can identify $\gnp$ with the product space $\Lambda_2\times\cdots\times\Lambda_n.$
Let
 $\pi_i:\gnp=\Lambda_2\times\cdots\times\Lambda_n\rightarrow\Lambda_i$
be the $i$th projection.
Let
 $A=\{G\in \gnp|\ \thet(G)\leq m\}$
and $B=\{H\in \gnp|\ m+\xi\leq\thet(H)\leq\thet_0\}$.

Let $H\in B$,
and let $(b_1,\ldots,b_n)$ be an orthogonal labeling of $\bar H$ such that
	\begin{equation}\label{eqLargedeva}
	m+\xi\leq\thet(H)=\sum_{i=1}^nc(b_i).
	\end{equation}
Set $\alpha_i=c(b_i)$, and $\alpha=(\alpha_2,\ldots,\alpha_n)$.
As $0\leq\alpha_i\leq1$ for all $i$, we have
 \begin{equation}\label{eqLargedevc}
 \sum_{i=2}^n \alpha_i^2\leq\sum_{i=1}^n\alpha_i=\thet(H)\leq\thet_0.
 \end{equation}
Now let $G\in A$, set $a_1=0$, and let
 $$a_i=\left\{\begin{array}{cl}b_i&\textrm{if }\pi_i(G)=\pi_i(H)\\0&\textrm{otherwise}
 \end{array}\right.\qquad\textrm{for $i=2,\ldots,n$.}$$
We claim that $(a_1,\ldots,a_n)$ is an orthogonal labeling of $\bar G$.
For if $i,j\in V$ are adjacent in $G$, and $i<j$, then we either have $\pi_j(G)=\pi_j(H)$
or $a_j=0$.
In the first case, $i$ and $j$ are adjacent in $H$, whence
 $\scal{a_i}{a_j}=\scal{b_i}{b_j}=0.$
Moreover, if $a_j=0$, then obviously $a_i\perp a_j$.
Thus, as $(a_1,\ldots,a_n)$ is an orthogonal labeling of $\bar G$, we have
	$\sum_{i=1}^nc(a_i)\leq\thet(G)\leq m.$
Hence, Eq.\ (\ref{eqLargedeva}) yields
 \begin{equation}\label{eqLargedevb}
 \xi\leq c(b_1)+\sum_{i=2}^nc(b_i)-c(a_i)\leq1+\sum_{i:\,\pi_i(G)\not=\pi_i(H)}c(b_i)
 	=1+\sum_{i:\,\pi_i(G)\not=\pi_i(H)}\alpha_i.
  \end{equation}

Set
 $t=(\xi-1)/\sqrt{\thet_0}.$
Then, by (\ref{eqLargedevc}) and (\ref{eqLargedevb}), for all $G\in A$ we have
 $$\sum_{i:\,\pi_i(G)\not=\pi_i(H)}\alpha_i\geq t\bc{\sum_{i=1}^n \alpha_i^2}^{1/2}.$$
Consequently, Thm.\ \ref{Thm_Talagrand} entails
 \begin{eqnarray*}
 \pr(A)\pr(m+\xi\leq\thet(\gnp)\leq\thet_0)&\leq&\exp(-t^2/4)
        =\exp\bc{-\frac{(\xi-1)^2}{4\thet_0}}\leq\exp\bc{-\frac{\xi^2}{5\thet_0}}.
 \end{eqnarray*}
Hence, our assertion follows from the fact that $\pr(A)\geq 1/2$.
\qed\end{proof}

\noindent\emph{Proof of Thm.\ \ref{Thm_Largedev}.}
By our assumption that $p\leq 0.99$, we have that $\thet(\gnp)\geq\alpha(\gnp)=\Omega(\ln(n))$ whp.
Hence, we can choose $n_0$ large enough such that any median $m$ of $\thet(\gnp)$ satisfies
$m\geq c_0$ for some large constant $c_0$.
As for the upper tail bound i.), by L.\ \ref{Lemma_Key} we have
	\begin{eqnarray*}
	\pr(m+\xi\leq\thet)&\leq&\sum_{l=1}^{\infty}\pr(m+l\xi\leq\thet\leq m+(l+1)\xi)\\
	&\leq&\sum_{l=1}^{\infty}2\exp\bc{-\frac{l^2\xi^2}{5(m+(l+1)\xi)}}
	\leq 2\sum_{l=1}^{\infty}\exp\bc{-\frac{l\xi^2}{5(m+2\xi)}}\\
	&\leq&30\exp\bc{-\frac{\xi^2}{5m+10\xi}},
\end{eqnarray*}
as desired.

In order to prove the lower tail bound ii.),
observe that by i.) we can choose $c_0$ large enough such that
$\pr(\thet(\gnp)\geq2m)<1/6$, say.
Consequently,
 $\pr(m\leq\thet(\gnp)\leq 2m)\geq 1/3.$
Let
 $A=\{G|\ \thet(G)\leq m-\xi\}$
and
 $B=\{H|\ m\leq\thet(H)\leq2m\}.$
Then, a similar argument as in the proof of L.\ \ref{Lemma_Key} yields
 $$\pr(A)\pr(B)\leq\exp(-t^2/4)=\exp\bc{-\frac{(\xi-1)^2}{8m}}\leq\exp\bc{-\frac{\xi^2}{10m}}.$$
Thus, our assertion follows from the fact that $\pr(B)\geq 1/3$.
\qed

\subsubsection{The large deviation result for $\thetI$.}
To prove the bounds in Thm.\ \ref{Thm_Largedev}
for $\thetI$, we make use of a characterization of $\thetI$ established in \cite{Groepl}.
Let $(v_1,\ldots,v_n)$ be an assignment of vectors $v_1,\ldots,v_n\in\RR^d$ to
the vertices of $G=(V,E)$, where $d>0$ is any integer.
Let us call $(v_1,\ldots,v_n)$ a \emph{strong orthogonal labeling} of $\bar G$ if
$v_i\perp v_j$ whenever $\{i,j\}\in E$, and
$\scal{v_i}{v_j}\geq0$ for all $i,j$.
As shown in \cite[pp.\ 51ff]{Groepl},
	\begin{equation}\label{eqClemens}
	\thetI(G)=\max\left\{\sum_{i=1}^nc(v_i)|\ (v_1,\ldots,v_n)\textrm{ is a strong orthogonal
		labeling of $\bar G$}\right\};
	\end{equation}
the proof goes along the lines of \cite{Knuth}.
Using Eq.\ (\ref{eqClemens}), the argument given for $\thet$ above carries
over without essential changes and yields the proof of the tail bounds
for $\thetI$.

\subsubsection{The large deviation result for $\thetII$.}
We shall establish a characterization of $\thetII$ that corresponds to
the characterization (\ref{eqorth}) of $\thet$, and which may be of independent interest.
Let $G=(V,E)$ be a graph.
If $x,y\in\RR^d$, then we let $c(x,y)=\scal{x}{y}^2\|x\|^{-2}\|y\|^{-2}$, if $x,y\not=0$,
and $c(x,y)=0$ otherwise.
Moreover, we call a family $(v_0,\ldots,v_n)$ of vectors $v_i\in\RR^d$ 
a \emph{weak orthogonal labeling} of $\bar G$ if
$\scal{v_0}{v_i}\geq 0$ for all $i$, and $\scal{v_i}{v_j}\leq0$ if $\{i,j\}\in E$, $i,j=1,\ldots,n$.
Here $d$ is any positive integer.
Note that a weak orthogonal labeling consists of $n+1=\#V+1$ vectors.
We define
	\begin{equation}\label{eqthetIIstrich}
	\thet_2'(G)=\max\left\{\sum_{i=1}^nc(v_i,v_0)|\ (v_0,\ldots,v_n)\textrm{ is a weak orthogonal labeling of }G\right\}.
	\end{equation}

\begin{lemma}
We have $\thetII(G)=\thet_2'(G)$ for all graphs $G=(V,E)$.
\end{lemma}
\begin{proof}
The following formulation of $\thetII$ as a semidefinite program
has been given in \cite{Szegedy}:
	\begin{equation}\label{eqthetIISDP}
	\thetII(G)=\max\sum_{i,j=1}^n b_{ij}\textrm{ s.t. }
		b_{ij}\leq0\textrm{ for all }\{i,j\}\in E\textrm{ and }\sum_{i=1}^nb_{ii}=1,
	\end{equation}
where the $\max$ is taken over all positive semidefinite matrices $B=(b_{ij})_{i,j}$.
(One can prove Eq.\ (\ref{eqthetIISDP}) e.g.\ using a similar argument as given in \cite{KMS}
to prove that (\ref{eqstrictvectorcol}) is equivalent to \Lovasz's original definition
of $\bthet$.)
To prove that $\thetII(G)\leq\thet_2'(G)$, let $B=(b_{ij})$ be a feasible matrix that
maximizes (\ref{eqthetIISDP}).
Since $B$ is positive semidefinite, there are vectors $b_1,\ldots,b_n\in\RR^n$ such that
$b_{ij}=\scal{b_i}{b_j}$.
Let $b=\sum_{i=1}^nb_i$.
Then $\scal{b_k}{b}\geq0$ for all $k$.
For assume otherwise, and
consider the matrix $B'=(b_{ij}')$, where $b_{ij}'=b_{ij}$ for $i,j\not=k$,
$b_{ki}'=b_{ik}'=0$ for $i\not=k$, and $b_{kk}'=b_{kk}$.
Then $B'$ is positive semidefinite, and is a feasible solution to (\ref{eqthetIISDP}).
Consequently, our assumption $0>\scal{b_k}{b}=\sum_jb_{kj}=\sum_jb_{jk}$ entails
	$$\sum_{i,j}b_{ij}=\thetII(G)\geq\sum_{i,j}b_{ij}'=\sum_{i,j}b_{ij}-2\sum_{i\not=k}b_{ik}
		>\sum_{i,j}b_{ij},$$
a contradiction.
Hence, $\scal{b_k}{b}\geq0$ for all $k$.
Letting $v_0=b/\|b\|$, $v_i=b_i/\|b_i\|$ if $b_i\not=0$, and $v_i=0$ otherwise ($i=1,\ldots,n$),
we obtain a weak orthogonal labeling of $\bar G$ satisfying
$\sum_{i=1}^nc(v_i,v_0)\geq\sum_{i,j}b_{ij}=\thetII(G)$ (cf.\ the proof of Thm.\ 5 in \cite{Lovasz}).

Conversely, let $(v_0,\ldots,v_n)$ be a weak orthogonal labeling of $\bar G$
such that $\sum_{i=1}^nc(v_i,v_0)=\thet_2'(G)$.
We may assume that $v_i$ either is a unit vector or is equal to zero for all $i$.
Let
	$$b_i=\thet_2'(G)^{-1/2}\scal{v_0}{v_i}v_i,$$
set $b_{ij}=\scal{b_i}{b_j}$, and $B=(b_{ij})_{i,j}$.
Then $B$ is positive semidefinite, and if $\{i,j\}\in E$, then
	$b_{ij}=\scal{v_0}{v_i}\scal{v_0}{v_j}\scal{v_i}{v_j}/\thet_2'(G)\leq0$,
because $\scal{v_0}{v_i},\scal{v_0}{v_j}\geq0\geq\scal{v_i}{v_j}$.
Moreover, $\sum_{i=1}^nb_{ii}=\sum_{i=1}^nc(v_0,v_i)/\thet_2'(G)=1$,
whence $B$ is a feasible solution to (\ref{eqthetIISDP}).
Finally, to show that $\sum_{i,j}b_{ij}\geq\thet_2'(G)$, we adapt the argument used in \cite{Groepl}
to prove (\ref{eqClemens}):
Let $M=\sum_{i=1}^nv_i{^tv_i}$.
Then
	$\thet_2'(G)=\sum_i\scal{v_0}{v_i}^2=\scal{Mv_0}{v_0}\leq\|Mv_0\|.$
Consequently,
	$$\thet_2'(G)\leq\|M(\thet_2'(G)^{-1/2}v_0)\|^2=
		\|\sum_{i=1}^n\thet_2'(G)^{-1/2}\scal{v_0}{v_i}v_i\|^2=\sum_{i,j}b_{ij}\leq\thetII(G),$$
thereby proving the lemma.
\qed\end{proof}

Using the characterization (\ref{eqthetIIstrich}) of $\thetII$, the arguments used
to prove Thm.\ \ref{Thm_Largedev} for $\thet(\gnp)$ also apply to $\thetII(\gnp)$.

\subsection{Concentration of $\vchi$, $\bthet$, and $\bthetII$ in Intervals
	of Constant Length}\label{Sec_Sharp}

Though the proofs of Thms.\ \ref{Thm_ShamirSpencer} and \ref{Thm_Luczak}
go along the lines of \cite{Luczak,ShamirSpencer},
we have to replace arguments concerning the chromatic
number by arguments that apply to $\vchi$, $\bthet$, and $\bthetII$.
We shall demonstrate the proofs for $\bthetII$, as this turns out
to be the most demanding case.
All arguments carry over to $\vchi$ and $\bthet$ immediately.
We adapt a simplification of the argument given in \cite{ShamirSpencer}
attributed to Frieze in \cite{Luczak}.

\subsubsection{Proof of Thm.\ \ref{Thm_ShamirSpencer}.}
Let $p$ and $\beta$ be as in Thm.\ \ref{Thm_ShamirSpencer}.
The proof is based on the following large deviation result,
which is a consequence of Azuma's inequality (cf.\ \cite[p.\ 37]{JLR}).

\begin{lemma}\label{Lemma_Azuma}
Suppose that $X:\gnp\rightarrow\RR$ is a random variable that satisfies
the following conditions for all graphs $G=(V,E)$.
\begin{itemize}
\item	For all $v\in V$ the following holds.
	Let $G^*=G+\{\{v,w\}|\ w\in V,\ w<v\}$, and let $G_*=G-\{\{v,w\}|\ w\in V,\ w<v\}$.
	Then $|X(G^*)-X(G_*)|\leq1$.
\item	If $H$ is a weak subgraph of $G$, then $X(H)\leq X(G)$.
\end{itemize}
Then $\pr(|X-\Erw(X)|>t\sqrt{n})\leq2\exp(-t^2/2)$.
\end{lemma}

Let $\omega=\omega(n)$ be a sequence tending to infinity slowly,
e.g.\ $\omega(n)=\ln\ln(n)$.
Furthermore, let
 \begin{equation}\label{eqk}
 k=k(n,p)=\inf\{x>0|\ \pr(\bthetII(\gnp)\leq x)\geq\omega^{-1}\}.
 \end{equation}
For any graph $G=(V,E)$ let
 $$Y(G)=Y_k(G)=\min\{\#U|\ U\subset V,\ \bthetII(G-U)\leq k\}.$$
Then $\bthetII(G)\leq k$ if and only if $Y(G)=0$.
Hence, $\pr(Y=0)\geq\omega^{-1}$.
Moreover, by Prop.\ \ref{Prop_facts} and (\ref{eqmonotone}), the random variable $Y$ satisfies the assumptions of L.\ \ref{Lemma_Azuma}.
Therefore, letting $\mu=\Erw(Y)$, for any $\lambda>0$ we have
	\begin{equation}\label{eqY}
	\pr(|Y(\gnp)-\mu|\geq\lambda\sqrt{n})\leq2\exp(-\lambda^2/2).
	\end{equation}
We claim that $\mu\leq\sqrt{n}\omega$.
For if $\mu>\sqrt{n}\omega$, then (\ref{eqY}) yields
 $$\omega^{-1}\leq\pr(Y=0)\leq\pr(Y\leq\mu-\sqrt{n}\omega)\leq2\exp(-\omega^2/2),$$
a contradiction.
Thus, again by L.\ \ref{Lemma_Azuma}, $Y\leq 2\sqrt{n}\omega$ with high probability.
The following lemma is implicit in \cite{ShamirSpencer} (cf.\ the proof of L.\ 8 in \cite{ShamirSpencer}).

\begin{lemma}\label{Lemma_Core}
Let $\delta>0$.
Whp.\ the random graph $G=\gnp$ enjoys the following property.
If $U\subset V$, $\#U\leq2\sqrt{n}\omega$,
then $\#E(G\lbrack U\rbrack)<\#Us/2$,
where $s>\frac{2}{2\beta-1}+\delta$.
Consequently, $\chi(G\lbrack U\rbrack)\leq s$,
\end{lemma}

To conclude the proof of Thm.\ \ref{Thm_ShamirSpencer},
let $G=\gnp$, and suppose that there is some $U\subset V$, $\#U\leq2\sqrt{n}\omega$,
such that $\bthetII(G-U)\leq k\leq\bthetII(G)$.
Since by L.\ \ref{Lemma_Core} $\bthetII(G\lbrack U\rbrack)\leq\chi(G\lbrack U\rbrack)\leq s$ whp.,
Prop.\ \ref{Prop_facts} entails that $k\leq\bthetII(G)\leq k+s$ whp.,
thereby proving Thm.\ \ref{Thm_ShamirSpencer}.

\subsubsection{Proof of Thm.\ \ref{Thm_Luczak}.}

Let $\omega=\omega(n)=(\ln\ln n)^{1/3}$ be a sequence tending to infinity slowly.
By L.\ \ref{Lemma_Core}, the random graph $G=\gnp$ admits
no $U\subset V$, $\#U\leq\omega^3\sqrt{n}$, spanning more than $3(\#U-\eps)/2$ edges whp.,
where $\eps>0$ is a small constant.
Let $k$ be defined as in (\ref{eqk}).
As shown in the proof of Thm.\ \ref{Thm_ShamirSpencer},
whp.\ there is a set $U\subset V$, $\#U\leq\omega\sqrt{n}$, such that $\bthetII(G-U)\leq k$.
Following \L uczak \cite{Luczak}, we let $U=U_0$, and construct a sequence $U_0,\ldots,U_m$ as follows.
If there is no edge $\{v,w\}\in E$ with $v,w\in N(U_i)\setminus U_i$, then we let $m=i$ and finish.
Otherwise, we let $U_{i+1}=U_i\cup\{v,w\}$ and continue.
Then $m\leq m_0=\omega^2\sqrt{n}$, because otherwise $\#U_{m_0}=(2+o(1))\omega^2\sqrt{n}$
and $\#E(G\lbrack U_{m_0}\rbrack)\geq 3(1-o(1))\#U_{m_0}/2$.
Let $R=U_m$.

By L.\ \ref{Lemma_Core}, $\bthetII(G\lbrack R\rbrack)\leq\chi(G\lbrack R\rbrack)\leq 3.$
Furthermore, $I=N(R)\setminus R$ is an independent set.
Let $G_1=G\lbrack R\cup I\rbrack$, $S=V\setminus(R\cup I)$, and $G_2=G\lbrack S\cup I\rbrack$.
Then $\bthetII(G_2)\leq k$, and $\bthetII(G_1)\leq 4$.
In order to prove that $\bthetII(G)\leq k+1$, we shall first construct a rigid vector
$k+1$-coloring of $G_2$ that assigns the same vector to all vertices in $I$.
Thus, let $(x_v)_{v\in S\cup I}$ be a rigid vector $k$-coloring of $G_2$.
Let $x$ be a unit vector perpendicular to $x_v$ for all $v\in S$.
Moreover, let $\alpha=(k^2-1)^{-1/2}$, and set
	$$y_v=\left\{\begin{array}{cl}(\alpha^2+1)^{-1/2}(x_v-\alpha x)&\textrm{ for $v\in S$}\\
			x&\textrm{ for }v\in I.\end{array}\right.$$
Then all $y_v$ are unit vectors, and if $v\in S$, $w\in I$, then
	$\scal{y_v}{y_w}=\scal{y_v}{x}=-1/k$.
Further, if $v,w\in S$ are adjacent in $G_2$, then
	$$\scal{y_v}{y_w}=\frac{1}{\alpha^2+1}\bc{\scal{x_v}{x_w}+\alpha^2}=
		\frac{1}{\alpha^2+1}\bc{-\frac{1}{k-1}+\alpha^2}=-\frac{1}{k}.$$
Likewise, if $v,w\in S$ are non-adjacent in $G_2$, then
$\scal{y_v}{y_w}\geq-1/k$, thereby proving that $(y_v)_{v\in S\cup I}$ is a
rigid vector $(k+1)$-coloring of $G_2$.
In a similar manner, we can construct a rigid vector $4$-coloring $(y_v')_{v\in R\cup I}$
of $G_1$ that assigns the same vector $x'$ to all vertices in $I$.

Applying a suitable orthogonal transformation if necessary, we may assume that $x=x'$.
Let $l=\max\{4,k+1\}$.
Since $N(R)\subset R\cup I$, we obtain a rigid vector $l$-coloring $(z_v)_{v\in V}$ of $G$,
where $z_v=y_v$ if $v\in S\cup I$, and $z_v=y_v'$ if $v\in R$.
By the lower bound on $\bthetII(\gnp)$ in Thm.\ \ref{Thm_thetgnp}
(which does not rely on Thm.\ \ref{Thm_Luczak} of course), choosing $c_0$
large enough we may assume that $k\geq 4$, whence $k\leq\bthetII(G)\leq k+1$.

\section{The Probable Value of $\thet(\gnp)$, $\bthet(\gnp)$, etc.}\label{Sec_thetgnp}

In Sec.\ \ref{Sec_Lower} we prove the lower bounds asserted in Thm.\ \ref{Thm_thetgnp}.
These follow from results on the SDP relaxation of MAX CUT on random graphs
due to Coja-Oghlan, Moore, and Sanwalani~\cite{MAXCUT}, and
do not depend on the concentration results in the previous section.
In Sec.\ \ref{Sec_UpperI} and \ref{Sec_UpperII} we prove the upper bounds on
$\thetII$ and $\bthetII$, which rely on Thm.\ \ref{Thm_Largedev}, Thm.\ \ref{Thm_Luczak},
and on a lemma on the spectrum of a certain auxiliary matrix given in Sec.\ \ref{Sec_Spec}.
Finally, in Sec.\ \ref{Sec_SDPk} we prove Cor.~\ref{Cor_SDPk}.

\subsection{The Lower Bound on $\vchi(\gnp)$}\label{Sec_Lower}

To bound $\bthetI(\gnp)$ from below, we make use of an estimate on the probable value
of the SDP relaxation $\SMC=\SDP_2$ of MAX CUT (cf.\ Sec.\ \ref{Sec_Pre} for the definition).
Combining Thms.\ 4 and~5 of \cite{MAXCUT} instantly yields the following bound on the probable value
of $\SMC(\gnp)$.

\begin{lemma}\label{Lemma_SMC}
Suppose that $c_0/n\leq p\leq1-c_0/n$ for some large constant $c_0>0$.
There is a constant $\lambda>0$ (independent of $n,p$) such that
	\begin{equation}\label{eqSMC}
	\pr\bc{\SMC(\gnp)>\frac{1}{2}\bink{n}{2}p+\lambda n^{3/2}p^{1/2}(1-p)^{1/2}}\leq\exp(-2n).
	\end{equation}
\end{lemma}

Let $G=(V,E)$ be a graph with adjacency matrix $A=(a_{ij})_{i,j=1,\ldots,n}$.
Let $v_1,\ldots,v_n$ be a vector $k$-coloring of $G$, where $k=\vchi(G)\geq2$.
Then $\|v_i\|=1$ for all $i$, and $\scal{v_i}{v_j}\leq-1/(k-1)$ whenever $\{i,j\}\in E$.
Therefore, we can interpret $v_1,\ldots,v_n$ as a feasible solution to $\SMC$, whence
	\begin{equation}\label{eqSMCvchi}
	\SMC(G)\geq\sum_{i<j}\frac{a_{ij}}{2}(1-\scal{v_i}{v_j})\geq\#E\bc{\frac{1}{2}+\frac{1}{k-1}}=
		\#E\bc{\frac{1}{2}+\frac{1}{\vchi(G)-1}}.
	\end{equation}
Let $c_0/n\leq p\leq 1-c_0/n$ for some large constant $c_0>0$.
By Chernoff bounds (cf.\ \cite[p.\ 26]{JLR}),
	\begin{equation}\label{eqConcE}
	\pr\bc{\#E(\gnp)<\bink{n}{2}p-8n^{3/2}p^{1/2}(1-p)^{1/2}}\leq\exp(-2n).
	\end{equation}
Combining (\ref{eqSMC}), (\ref{eqSMCvchi}), and (\ref{eqConcE}), we conclude that
	$$\vchi(\gnp)\geq\vchi(\gnp)-1\geq
		\frac{\bink{n}{2}p-8n^{3/2}p^{1/2}(1-p)^{1/2}}{(\lambda+4)n^{3/2}p^{1/2}(1-p)^{1/2}}
			\geq\frac{1}{2(\lambda+4)}\sqrt{\frac{np}{1-p}}$$
holds with probability at least $1-\exp(-n)$.
As $\bar\gnp=G_{n,1-p}$, this proves (\ref{eqexp}) and the lower bounds in Thm.\ \ref{Thm_thetgnp}.

\begin{remark}
Suppose that $np=O(1)$.
Then (\ref{eqexp}) shows that the probability that $\vchi(\gnp)$ is less than
$c\sqrt{np}$ is exponentially small, for some constant $c>0$.
It is easily seen that no similar statement holds for the upper tail,
i.e.\ for the event that $\vchi(\gnp)>\zeta\sqrt{np}$ for some large $\zeta$.
The reason is that the probability that $\omega(\gnp)>\zeta\sqrt{np}$ is at least
$p^{\zeta^2np}\geq\exp(-O(\ln(n))$, and $\vchi(G)\geq\omega(G)$ for all $G$.
\end{remark}

\subsection{Spectral Considerations}\label{Sec_Spec}

Let us briefly recall \Juhasz's proof that $\thet(\gnp)\leq(2+o(1))\sqrt{n(1-p)/p}$
for constant values of $p$, say.
Given a graph $G=(V,E)$, we consider the matrix $M=M(G)=(m_{ij})_{i,j=1,\ldots,n}$, where
	\begin{equation}\label{eqdefMstrich}
	m_{ij}=\left\{\begin{array}{cl}1&\textrm{ if }\{i,j\}\not\in E\\(p-1)/p&\textrm{ otherwise,}
		\end{array}\right.\qquad(i\not=j),
	\end{equation}
and $m_{ii}=1$ for all $i$.
Then $\lambda_1(M)\geq\thet(G)$.
Moreover, as $p$ is constant,
the result of F\"uredi and Komlos \cite{Komlos} on the eigenvalues
of random matrices applies and yields that $\thet(\gnp)\leq\lambda_1(M)\leq(2+o(1))\sqrt{n(1-p)/p}$ whp.
This argument carries over to the case $\ln(n)^7/n\leq p\leq1/2$:

\begin{lemma}\label{Lemma_eigdense}
Suppose that $\ln(n)^7/n\leq p\leq1/2$.
Then $\|M(\gnp)\|\leq 3\sqrt{n/p}$ whp.
\end{lemma}
\begin{proof}
As it is assumed in \cite{Komlos} that the variance of the matrix entries is independent
of $n$, the proof of F\"uredi and Komlos \cite{Komlos} needs some minor adaptions
to prove the lemma;
all details have been carried out in \cite[Sec.\ 4]{CojaTaraz}.
\qed\end{proof}

However, it is easily seen that in the sparse case, e.g.\ if $np=O(1)$,
we have $\lambda_1(M)\gg n$ whp.
The reason is that in the case $np\geq\ln(n)^7$ the random graph $\gnp$
is ``almost regular'', which is not true if $np=O(1)$ (cf.\ \cite{KrivSudakov}).
We will get around this problem by chopping off all vertices of degree
considerably larger than $np$, as first proposed in \cite{AlonKahale}.
Thus, let $\eps>0$ be a small constant,
and consider the graph $G'=(V',E')$ obtained from $G=\gnp$
by deleting all vertices of degree greater than $(1+\eps)np$.

\begin{lemma}\label{Lemma_Mstrich}
Suppose that $c_0/n\leq p\leq\ln(n)^7/n$ for some large constant $c_0$.
Let $G=\gnp$, and let $M'=M(G')$.
Then $\pr(\|M'\|\leq c_1\sqrt{n/p})\geq9/10,$
where $c_1>0$ denotes some constant.
\end{lemma}

To prove L.\ \ref{Lemma_Mstrich}, we make use of
the following lemma, which is implicit in~\cite[Sections~2~and~3]{Ofek}.

\begin{lemma}\label{Lemma_Ofek}
Let $G=\gnp$ be a random graph, where $c_0/n\leq p\leq\ln(n)^7/n$ for some large constant $c_0>0$.
Let $n'=\#V(G')$, $e={n'}^{-1/2}\vecone\in\RR^{n'}$, and $A'=A(G')$.
For each $\delta>0$ there is a constant $C(\delta)>0$ such that in the case $np\geq C(\delta)$
with probability $\geq1-\delta$ we have
	\begin{equation}\label{eqOfek}
	\max\{|\scal{A'v}{e}|,|\scal{A'v}{w}|\}\leq c_1\sqrt{np}
		\textrm{ for all $v,w\perp\vecone$, $\|v\|=\|w\|=1$}.
	\end{equation}
Here $c_1>0$ denotes a certain constant.
\end{lemma}

In addition, the proof of Lemma \ref{Lemma_Mstrich} needs the following observation.

\begin{lemma}\label{Lemma_bard}
Let $c_1$ be a large constant.
The probability that in $G=\gnp$ there exists a set $U\subset V$, $\#U\geq n/2$, such that
$|\#E(G\lbrack U\rbrack)-\#U^2p/2|\geq c_1(\#U)^{3/2}p^{1/2}$ is less than $\exp(-n)$.
\end{lemma}
\begin{proof}
There are at most $2^n$ sets $U$.
By Chernoff bounds (cf.\ \cite[p.\ 26]{JLR}), for a fixed $U$ the probability that
$|\#E(G\lbrack U\rbrack)-\#U^2p/2|\geq c_1(\#U)^{3/2}p^{1/2}$
is at most $\exp(-2n)$, provided that $c_0$, $c_1$ are large enough.
\qed\end{proof}

%\medskip%
\noindent\emph{Proof of Lemma \ref{Lemma_Mstrich}.}
Let $G=\gnp$, let $n'=\#V(G')$, and let $A'$, $e$ be as in L.\ \ref{Lemma_Ofek}.
Without loss of generality, we may assume that $V'=V(G')=\{1,\ldots,n'\}$.
Let $c_1>0$ be a sufficiently large constant.
Let $J$ signify the $n'\times n'$ matrix with all entries equal to $1$.
Letting $\delta>0$ be sufficiently small and $c_0\geq C(\delta)$,
we assume in the sequel that (\ref{eqOfek}) holds, and that
$G$ has the property stated in L.\ \ref{Lemma_bard}.
Let $z\in\RR^{n'}$, $\|z\|=1$.
Then we have a decomposition $z=\alpha e+\beta v$, $\|v\|=1$, $v\perp\vecone$,
$\alpha^2+\beta^2=1$.
Since $\|M'z\|\leq\|M'e\|+\|M'v\|$,
if suffices bound $\max_{v\perp e,\|v\|=1}\|M'v\|$ and $\|M'e\|$.

Let $\rho:\RR^{n'}\rightarrow\RR^{n'}$ be the projection on the space $\vecone^{\perp}$.
Then
	$A'v=\rho A'v+\scal{A'v}{e}e,$
whence
	$\|A'v\|\leq\|\rho A'v\|+c_1\sqrt{np},$
for all unit vectors $v\perp\vecone$.
In order to bound $\|\rho A'v\|$, we estimate $\|\rho A'\rho\|$ via (\ref{eqOfek}):
	$$\|\rho A'\rho\|=\sup_{\|y\|=1}|\scal{\rho A'\rho y}{y}|
		=\sup_{\|y\|=1}|\scal{A'\rho y}{\rho y}|=\sup_{\|y\|=1,\ \vecone\perp y}|\scal{A'y}{y}|
			\leq c_1\sqrt{np}.$$
Consequently,
	$\|M'v\|=\|(J-\frac{1}{p}A')v\|=\frac{1}{p}\|A'v\|\leq 2c_1\sqrt{n/p}$
for all unit vectors $v\perp\vecone$.

To bound $\|M'e\|$, note that $-pM'=A'-pJ$.
Let $\bar d=2\#E(G')/n'$, and $x=A'e-(\bar d/n')Je$.
Then $x\perp\vecone$, and by (\ref{eqOfek}) we have
	$\|x\|^2=\scal{A'e}{x}-\scal{(\bar d/n')Je}{x}=
		\scal{A'e}{x}\leq c_1\sqrt{np}\|x\|,$
whence $\|x\|\leq c_1\sqrt{np}$.
By L.\ \ref{Lemma_bard}, $|\bar d-n'p|\leq c_1\sqrt{np}$.
As a consequence,
	$\|(\bar d/n')Je-pJe\|\leq c_1\sqrt{np}.$
Therefore,
	$\|pM'e\|\leq\|x\|+\|(\bar d/n')Je-pJe\|\leq 2c_1\sqrt{np},$
i.e.\ $\|M'e\|\leq2c_1\sqrt{n/p}$.
\qed

\subsection{Bounding $\thetII(\gnp)$ from above}\label{Sec_UpperI}

\emph{Let $c_0/n\leq p\leq1/2$ for some large constant $c_0>0$.}
First we observe that the largest eigenvalue of the matrix $M(G)$ considered
in the previous section provides an upper bound on $\thetII(G)$.
(Actually this follows from the characterization of $\bthetII$ as an eigenvalue
minimization problem given in \cite{Szegedy}.
However, as \cite{Szegedy} does not contain the proof, we show a brief ad hoc argument.)

\begin{lemma}\label{Lemma_eigbound}
Let $G$ be any graph.
Let $M=M(G)$.
Then $\lambda_1(M)\geq\bthetII(G)$.
\end{lemma}
\begin{proof}
Let $\lambda>\lambda_1(M)$.
Then the matrix $\lambda E_n-M$ is positive definite,
whence there exist vectors $b_1,\ldots,b_n\in\RR^n$ such
that $m_{ij}=-\scal{b_i}{b_j}$ for $i\not=j$, and
	$\lambda-1=\lambda-m_{ii}=\scal{b_i}{b_i}=\|b_i\|^2>0.$
Let $a_i=(\lambda-1)^{-1/2}b_i$.
Then $\|a_i\|=1$ for all $i$.
Moreover, if $i\not=j$ and $\{i,j\}\not\in E$, then
	$\scal{a_i}{a_j}=m_{ij}/(\lambda-1)=-1/(\lambda-1)$.
If $\{i,j\}\in E$, then $\scal{a_i}{a_j}\geq0$.
Hence $(a_1,\ldots,a_n)$ is a rigid vector $\lambda$-coloring of $G$.
Therefore, $\thetII(\gnp)\leq\lambda$ for all $\lambda>\lambda_1(M)$.
\qed\end{proof}

In the case $\ln(n)^7/n\leq p\leq1/2$, combining L.\ \ref{Lemma_eigdense} and L.~\ref{Lemma_eigbound}
yields that $\thetII(\gnp)\leq c_2\sqrt{n/p}$ whp.\ for some constant $c_2>0$, as desired.
Thus, let us assume that $c_0/n\leq p\leq\ln(n)^7/n$ in the sequel.
Let $\eps>0$ be a small constant.

\begin{lemma}\label{Lemma_highdegree}
With probability at least $9/10$ the random $\gnp$  has at most $1/p$ vertices
of degree greater than $(1+\eps)np$.
\end{lemma}
\begin{proof}
For each vertex $v$ of $\gnp$, the degree $d(v)$ is binomially distributed with
mean $(n-1)p$.
By Chernoff bounds (cf.\ \cite[p.\ 26]{JLR}), the probability that $d(v)>(1+\eps)np$ is at most
$\exp(-\eps^2np/100)$.
Hence, the expected number of vertices $v$ such that
$d(v)>(1+\eps)np$ is at most
	$n\exp(-\eps^2np/100)<1/(10p),$
provided $np\geq c_0$ for some large constant $c_0>0$.
Therefore, the assertion follows from Markov's inequality.
\qed\end{proof}

Let $G=\gnp$, and let $G'=(V',E')$ be the graph obtained from $G$ by deleting
all vertices of degree greater than $(1+\eps)np$.
Let $V''=V\setminus V'$, and $G''=G\lbrack V''\rbrack$.
Combining L.\ \ref{Lemma_highdegree} and L.~\ref{Lemma_Mstrich},
we obtain that 
	$\pr\bc{\thetII(G')\leq c_2\sqrt{n/p}\textrm{ and }
		\thetII(G'')\leq\#V(G'')\leq1/p\leq\sqrt{n/p}}>1/2,$
where $c_2$ denotes a suitable constant.
Consequently, Prop.\ \ref{Prop_facts} yields that
	$$\pr(\thetII(\gnp)\leq(c_2+1)\sqrt{n/p})>1/2.$$
Let $\mu=(c_2+1)\sqrt{n/p}$, $t=\ln(n)\sqrt{n}$, and note that $t=o(\sqrt{n/p})$.
Then, by Thm.\ \ref{Thm_Largedev},
	$$\pr(\thetII(\gnp)>\mu+t)\leq 30\exp\bc{-\Omega(\ln(n)^2)}=o(1).$$
Since $t<\sqrt{n/p}$, we get that $\thetII(\gnp)\leq(c_2+2)\sqrt{n/p}$ with high probability.

\subsection{Bounding $\bthetII(\gnp)$ from Above}\label{Sec_UpperII}

Let us first assume that $\ln(n)^7/n\leq p\leq1/2$.
Let $G=(V,E)=\gnp$ be a random graph, and consider the matrix
$\bar M=\frac{1}{1-p}E_{n}-\frac{p}{1-p}M(G),$
where $E_n$ is the $n\times n$-unit matrix,
and $M(G)$ is the matrix defined in (\ref{eqdefMstrich}).
Combining L.\ \ref{Lemma_eigdense} and L.\ \ref{Lemma_eigbound}, we have
	$$\bthetII(G)\leq\lambda_1(\bar M)\leq\|\frac{1}{1-p}E-\frac{p}{1-p}M'\|
		\leq\frac{p}{1-p}\|M\|+2\leq c_4\sqrt{np}$$
whp., where $c_4>0$ is a certain constant.

Now let $c_0/n\leq p\leq\ln(n)^7/n$ for some large constant $c_0>0$.
In this case, the proof of our upper bound on $\bthetII(\gnp)$
relies on the concentration result Thm.\ \ref{Thm_Luczak}.

\begin{lemma}\label{Lemma_bthetupper}
Whp.\ the random graph $G=\gnp$ admits no set $U\subset V$,
$\#U\leq1/p$, such that $\chi(G\lbrack U\rbrack)>\sqrt{np}$.
\end{lemma}
\begin{proof}
We shall prove that for all $U\subset V$, $\#U=\nu\leq1/p$, we have
$\#E(G\lbrack U\rbrack)<\nu\sqrt{np}/2$.
Then each subgraph $G\lbrack U\rbrack$ has a vertex of degree $<\sqrt{np}$,
a fact which immediately implies our assertion.
Thus, let $\nu\leq1/p$.
The probability that there exists some $U\subset V$, $\#U=\nu$,
$\#E(G\lbrack U\rbrack)\geq\nu\sqrt{np}/2$, is at most
	$$ %\begin{equation*}\label{eqCore}
	\bink{n}{\nu}\bink{\bink{\nu}{2}}{\nu\sqrt{np}/2}p^{\nu\sqrt{np}/2}
	\leq\bc{\frac{\eul n}{\nu}\bcfr{\eul\nu^2p}{\nu\sqrt{np}}^{\sqrt{np}/2}}^{\nu}
	=\bc{\frac{\eul n}{\nu}\bcfr{\eul\nu\sqrt{p}}{\sqrt{n}}^{\sqrt{np}/2}}^{\nu}
	$$%\end{equation*}
Let $b_{\nu}=(\eul n/\nu)(\eul\nu\sqrt{p}/\sqrt{n})^{\sqrt{np}/2}$.
Observe that the sequence $(b_{\nu})_{\nu=1,\ldots,n}$ is monotone increasing,
and that $b_{1/p}=\eul np(\eul/\sqrt{np})^{\sqrt{np}/2}\leq\exp(-2)$.
Therefore,
	$$\sum_{\nu=\ln(n)}^{1/p}b_{\nu}^{\nu}\leq b_{1/p}^{\ln(n)}/p\leq n^{-2}p^{-1}=o(1).$$
Moreover, if $\nu\leq\ln(n)$, then
	$b_{\nu}\leq\eul n\nu^{-1}(\eul\nu\sqrt{p}/\sqrt{n})^{\sqrt{np}/2}\leq1/n$,
whence $\sum_{\nu=1}^{\ln n}b_{\nu}^{\nu}=o(1)$.
Thus, $\sum_{\nu=1}^{1/p}b_{\nu}^{\nu}=o(1)$,
thereby proving the lemma.
\qed\end{proof}

Let $G=(V,E)=\gnp$ be a random graph, and let $G'=(V',E')$ be the graph obtained from $G$ by
removing all vertices of degree greater than $(1+\eps)np$, where $\eps>0$ is small
but constant.
Let $V''=V\setminus V'$, and let $G''=G\lbrack V''\rbrack$.
By L.\ \ref{Lemma_highdegree}, with probability at least $9/10$ we have $\#V''\leq1/p$.
Therefore, by L.\ \ref{Lemma_bthetupper},
	$$\pr(\bthetII(G'')\leq\sqrt{np})\geq\pr(\chi(G'')\leq\sqrt{np})\geq9/11.$$
To bound $\bthetII(G')$, we consider the matrix
	$\bar M=\frac{1}{1-p}E_{n'}-\frac{p}{1-p}M(G'),$
where $E_{n'}$ is the $\#V'\times\#V'$-unit matrix,
and $M(G')$ the matrix (\ref{eqdefMstrich}).
By L.\ \ref{Lemma_eigbound}, $\bthetII(G')\leq\lambda_1(\bar M)$.
Moreover, by L.\ \ref{Lemma_Mstrich},
with probability $\geq9/10$ we have
	$$\bthetII(G')\leq\lambda_1(\bar M)\leq\|\frac{1}{1-p}E-\frac{p}{1-p}M'\|
		\leq\frac{p}{1-p}\|M'\|+2\leq c_4\sqrt{np},$$
for some constant $c_4>0$.
Prop.\ \ref{Prop_facts} implies that $\bthetII(G)\leq\bthetII(G')+\bthetII(G'')$,
whence we conclude that
	$\pr(\bthetII(\gnp)\leq (c_4+1)\sqrt{np})>1/2.$
Since Thm.\ \ref{Thm_Luczak} shows that $\bthetII(\gnp)$ is concentrated in width one,
we have
	$$\pr\bc{\vchi(\gnp)\leq\bthet(\gnp)\leq\bthetII(\gnp)\leq (c_4+1)\sqrt{np}+1}=1-o(1),$$
thereby completing the proof of Thm.\ \ref{Thm_thetgnp}.

\begin{remark}
One could prove slightly weaker results on the probable value of $\thet(\gnp)$
and $\bthet(\gnp)$ than provided by Thm.\ \ref{Thm_thetgnp} without applying
any concentration results, or bounds on the SDP relaxation $\SMC$ of MAX CUT.
Indeed, using only Lemmata\ \ref{Lemma_highdegree}, \ref{Lemma_bthetupper}, \ref{Lemma_Mstrich}
(thus implicitly~\cite{Ofek}) and the estimates proposed in \cite{Juhasz},
one could show that for each $\delta>0$ there is $C(\delta)>0$ such the following holds.
If $np\geq C(\delta)$, then
	\begin{equation}\label{eqClaimOfek}
	\pr(c_1\sqrt{n/p}\leq\thet(\gnp)\leq c_2\sqrt{n/p})\geq1-\delta,\ \pr(c_3\sqrt{np}\leq\bthet(\gnp)\leq c_4\sqrt{np})\geq1-\delta.
	\end{equation}
Such an approach is mentioned without proof independently in the latest version of~\cite{Ofek}
(the phrase ``with high probability'' is used in the sense ``with probability $1-o(1)$ as $np\rightarrow\infty$''
in that paper).
However, Thm.\ \ref{Thm_thetgnp} is a bit stronger than (\ref{eqClaimOfek}),
as the bounds (\ref{eqthetgnp}) on $\thet(\gnp)$ and $\bthet(\gnp)$ hold with probability $1-o(1)$ as
$n\rightarrow\infty$ even if $np$ remains bounded.
Moreover, using the simpler
approach it seems hard to obtain exponentially small
probabilities as in (\ref{eqexp}).
\end{remark}

\subsection{The Lower Bound on $\SDP_k(\gnp)$}\label{Sec_SDPk}

Having established Thm.\ \ref{Thm_thetgnp}, we know that there exist
constants $c_0,c_1,c_2>0$ such that in the case $c_0/n\leq p\leq 1/2$ we have
	\begin{equation}\label{eqbthetIISDPk}
	c_1\sqrt{np}\leq\bthetII(\gnp)\leq c_2\sqrt{np}
	\end{equation}
with high probability.
Let $k\geq2$ be a fixed integer, and let us assume that $c_3k^2/n\leq p\leq1/2$,
where $c_3=\max\{c_0,c_1^{-1}\}$.
Let $G=\gnp$ satisfy (\ref{eqbthetIISDPk}), and consider an
rigid vector $\bthetII(G)$-coloring $(v_1,\ldots,v_n)$ of $G$.
Then 
	$$\scal{v_i}{v_j}\geq-\frac{1}{c_1\sqrt{np}-1}\geq-\frac{1}{k-1}$$
for all $i,j$, whence $(v_1,\ldots,v_n)$ is a feasible solution to $\SDP_k$.
Furthermore, if $\{i,j\}\in E$, then
	$$\scal{v_i}{v_j}\leq-\frac{1}{c_2\sqrt{np}-1}.$$
Consequently, letting $A=(a_{ij})_{i,j=1,\ldots,n}$ be the adjacency matrix of $G$,
we have
	\begin{eqnarray*}
	\SDP_k(G)&\geq&\sum_{i<j}a_{ij}\frac{k-1}{k}(1-\scal{v_i}{v_j})
		\geq\bc{1-\frac{1}{k}}\#E(G)+\frac{\#E(G)}{2c_2\sqrt{np}}.
	\end{eqnarray*}
As $\#E(\gnp)$ is concentrated about its mean $\bink{n}{2}p$,
we conclude that
	$$\SDP_k(\gnp)\geq\bc{1-\frac{1}{k}}\bink{n}{2}p+\frac{n^{3/2}p^{1/2}}{3c_2}$$
with high probability, thereby proving Cor.\ \ref{Cor_SDPk}.

\begin{remark}
Consider the following relaxation $\SDP_k'$ of $\SDP_k$:
	$$\SDP_k'(G)=\max\sum_{i<j}a_{ij}\frac{k-1}{k}\bc{1-\scal{v_i}{v_j}}\textrm{ s.t.\ }\|v_i\|=1,$$
where $A=(a_{ij})_{i,j=1,\ldots,n}$ is the adjacency matrix of $G$ and the $\max$ is taken
over all families $v_1,\ldots,v_n$ of unit vectors in $\RR^n$.
Then $\SDP_k'(G)=(2(k-1)/k)\SMC(G)$.
Consequently, L.~\ref{Lemma_SMC} shows that $\SDP_k'(\gnp)\leq(1-1/k)\bink{n}{2}p+c_3n^{3/2}p^{1/2}$ whp.,
where $c_3>0$ is some constant.
Thus, Cor.\ \ref{Cor_SDPk} implies that whp.\ both $\SDP_k$ and $\SDP_k'$ overestimate the weight
$\mkc(\gnp)$ of a \MKC\ by at least $c_4n^{3/2}p^{1/2}$, for some constant $c_4>0$.
Thus, in the case of random graphs the additional constraints $\scal{v_i}{v_j}\geq-1/(k-1)$ only affect
the precise constant in front of the second order term $n^{3/2}p^{1/2}$.
\end{remark}

\section{Random Regular Graphs}\label{Sec_reg}

We show how to adapt the arguments given in the previous section
to cover the case of random regular graphs.
Throughout we assume that $r\geq c_0$ for some large constant $c_0>0$.

The proof of (\ref{eqvchignr}) relies on the upper bound
on $\SMC(\gnr)$ \cite[Thm.~15]{MAXCUT}, and is similar to the proof of (\ref{eqexp}).
To prove the upper bound on $\thetII(\gnr)$, $c_0\leq r=o(n^{1/4})$, we switch to
the \emph{configuration model} (cf.\ \cite{Wormald}).
Let $W=V\times\{1,\ldots,r\}$.
The elements of $W$ are called \emph{half edges}.
A \emph{configuration} $\rho$ is a partition of $W$ into $m=rn/2$ pairs, where we assume that $rn$ is even.
Thus, to each half-edge $(u,v)$, $\rho$ assigns another half-edge $\rho(u,v)\not=(u,v)$
such that $\rho^2=\mathrm{id}$.
We say that $(u,v)$ and $\rho(u,v)$ form an \emph{edge}.
By $\conf=\conf(r)$ we denote the set of all configurations.
Then $\#\conf(r)=(2m-1)!!$.

To each $\rho\in\conf$,
the canonical map $\pi:W\rightarrow V$ assigns an $r$-regular multigraph $\pi(\rho)$.
If we equip $\conf$ with the uniform distribution,
then conditional on $\gnr$, $\pi$ induces the uniform distribution.
By $\sigma(\rho)$ we denote the simple graph obtained from the multigraph $\pi(\rho)$
by deleting all loops and turning all multiple edges into single edges.
We define the \emph{adjacency matrix} $A=A(\rho)=(a_{ij})_{i,j}$
of $\rho\in\conf$ to be the matrix with entries
	$$a_{ij}=\textrm{number of edges joining $i$ and $j$ in $\pi(\rho)$}$$
if $i\not=j$, and let $a_{ii}$ be twice the number of loops at vertex $i$ in $\pi(\rho)$.

\begin{lemma}\label{Lemma_ConcReg}
Let $i\in\{1/2,1,2\}$.
Let $\mu$ be the expectation of $\thet_i\circ\sigma$ over $\conf$
(where $\thet_i\circ\sigma(\rho)=\thet_i(\sigma(\rho))$).
Then for any $t>0$ we have
	$\pr(|\thet_i\circ\sigma-\mu| > t ) \leq 2\exp(-t^2/(128m))\enspace .$
\end{lemma}
\begin{proof}
If two configurations $\rho$, $\rho'$ are such that $\sigma(\rho)$ can be transformed
into $\sigma(\rho')$ by adding or deleting at most $4$ edges, then
$|\thet_i\circ\sigma(\rho)-\thet_i\circ\sigma(\rho')|\leq4$ by Prop.\ \ref{Prop_facts}.
Therefore, the martingale argument used in the proof of Lemma 14 of \cite{MAXCUT}
carries over without essential changes.
\qed\end{proof}

The proof of the following lemma goes along the lines of \cite{Kahn}.
However, as we work with the configuration model and consider also the case
that the degree $r$ tends to infinity, some adaptions are necessary;
these have been carried out in \cite{MAXCUT}.

\begin{lemma}\label{Lemma_Kahn}
There is a constant $\gamma>0$ such that
with high probability the adjacency matrix $A=A(\pi(\rho))$, $\rho\in\conf$,
satisfies $|\scal{Ax}{y}|\leq\gamma\sqrt{r}$,
for all unit vectors $x\perp\vecone$, $y\perp\vecone$.
\end{lemma}

Given a configuration $\rho\in\conf$, we let $M=M(\rho)=J-\frac{n}{r}A(\pi\rho)$.
Then $M$ is a symmetric $n\times n$ matrix.

\begin{lemma}\label{Lemma_MSpec}
There is a constant $c>0$ such that whp.\ we have
$\|M(\rho)\|\leq c nr^{-1/2}$.
\end{lemma}
\begin{proof}
Let $M=M(\rho)$.
As $M\vecone=0$, it suffices to prove that $\scal{Mv}{w}\leq cnr^{-1/2}$ for all
unit vectors $v,w\perp\vecone$ whp.
But this follows from L.\ \ref{Lemma_Kahn} easily.
\qed\end{proof}

\begin{lemma}\label{Lemma_Loop}
The expected number of loops and multiple edges in $\pi(\rho)$, $\rho\in\conf$,
is at most $9r^2$.
Hence, with probability $\geq1/2$ there are at most $18r^2$ loops or multiple edges.
\end{lemma}
\begin{proof}
Let $(u,t),(u,s),(v,r),(v,r')\in W$ be distinct half edges.
The probability that $\rho(u,t)=(v,r)$ and $\rho(u,s)=(v,r')$ is $\sim(2m)^{-2}$.
There are $2m$ choices of $(v,r)$,
and then at most $r$ choices of $(v,r')$.
Further, given $(u,t)$, there are $r$ possible choices of $s$.
Hence, the expected number of half edges that participate in multiple edges is
at most $2(2m)^2r^2(2m)^{-2}\leq8r^2$.

As for loops, let $(u,t)\in W$, and let $s\in\{1,\ldots,r\}\setminus\{t\}$.
The probability that $\rho(u,t)=(u,s)$ is $\sim(2m)^{-1}$.
Since there are at most $2m$ possible choices of $(u,t)$, and then at most $r$ choices of $s$,
the expected number of loops is at most $2(2m)r(2m)^{-1}\leq4r<r^2$.
\qed\end{proof}

\begin{lemma}\label{Lemma_almostreg}
There are constant $c_1,c_2>0$ such that whp.\
a random configuration $\rho\in\conf$ satisfies
$c_1 nr^{-1/2}\leq\thetI(\sigma\rho)\leq\thet(\sigma\rho)\leq\thetII(\sigma\rho)\leq c_2 nr^{-1/2}$.
\end{lemma}
\begin{proof}
Let $\mathcal{B}$ be the event that the number of multiple edges and loops in $\pi\rho$,
$\rho\in\conf$, is at most $20r^2$.
By L.\ \ref{Lemma_Loop}, $\pr(\mathcal{B})\geq1/2$.
Consequently, by L.\ \ref{Lemma_MSpec}, there is a constant $c_1>0$ such that
	$\pr(\|M(\rho)\|\leq c_1 nr^{-1/2}|\mathcal{B})\geq1/2.$
Hence, $\pr(\rho\in\mathcal{B}\textrm{ and }\|M(\rho)\|\leq c_1 nr^{-1/2})\geq 1/4$.

We claim that if $\rho\in\mathcal{B}$ satisfies $\|M(\rho)\|\leq c_1 nr^{-1/2}$,
then $\thetII(\sigma\rho)\leq 2c_1 nr^{-1/2}$.
For let $Y$ be the set of all vertices $v\in V$ that participate in a multiple edge or a loop.
Then $y=\#Y\leq 40r^2$.
Relabeling the vertices if necessary, we may assume that $Y=\{n-y+1,\ldots,n\}$.
Let $M=M(\rho)=(m_{ij})_{i,j=1,\ldots,n}$, and set $M'=(m_{ij})_{i,j=1,\ldots,n-y}$.
Then
	$$m_{ij}=\left\{\begin{array}{cl}1&\textrm{ if $i,j$ are non-adjacent in $\sigma\rho$}\\
		(r-n)/r&\textrm{ otherwise,}\end{array}\right.\qquad(1\leq i<j\leq n-y),$$
and $m_{ii}=1$ for all $i$.
Let $H$ be the simple graph on $V(H)=\{1,\ldots,n-y\}$ induced by $\sigma\rho$.
Then
	$\thetII(H)\leq\lambda_1(M')\leq\|M\|\leq c_1nr^{-1/2}.$
Since the graph $\sigma\rho$ can be obtained from $H$ by adding $y$ vertices,
and since $y\leq nr^{-1/2}$, we conclude that
	$\thetII(\sigma\rho)\leq\thetII(H)+y\leq  2c_1 nr^{-1/2}.$
Hence, $\pr(\thetII(\sigma\rho)\leq 2c_1 nr^{-1/2})\geq1/4$.
Invoking L.\ \ref{Lemma_ConcReg} completes the proof of the upper bound.

As for the lower bound, let $c_2$ be a sufficiently large constant, and let
$\rho\in\mathcal{B}$ be such that the adjacency matrix $A=A(\pi(\rho))$
satisfies $|\scal{A\xi}{\eta}|\leq c_2r^{1/2}$ for all unit vectors $\xi,\eta\perp\vecone$.
Let $Y$, $y$, and $H$ be as before, $y\leq40r^2$.
Moreover, let $\bar A$ be the adjacency matrix of $\overline{\sigma(\rho)}\lbrack V\setminus Y\rbrack$,
and let $E_{n-y}$ and $E_n$ be unit matrices of size $n-y$ and $n$.
Since $\vecone$ is an eigenvector of $J-A+E_n$, we have
	$$\lambda_n(J-A+E_n)=1-\lambda_2(A)=\max_{\xi\perp\vecone,\ \|\xi\|=1}\scal{A\xi}{\xi}\geq-c_2r^{1/2}.$$
Hence,
	$\lambda_n(\bar A)\geq-c_2r^{1/2},$
because $\bar A=J-A(\sigma(\rho)\lbrack V\setminus Y\rbrack)+E_{n-y}$ is a principal minor
of $J-A+E_n$.
Let $B=(n-y)^{-1}\bc{E_{n-y}-\lambda_n(\bar A)^{-1}\bar A}$.
Then the matrix $B=(b_{ij})_{i,j=1,\ldots,n-y}$ is positive semidefinite,
and we have $\sum_ib_{ii}=1$.
Moreover, $b_{ij}\geq0$ for all $i,j$, and
if $i,j\in V\setminus Y$ are adjacent in $\sigma(\rho)$, then $b_{ij}=0$.
It is shown in \cite[pp.\ 51ff]{Groepl} that such a matrix $B$ satisfies
$\thetI(\sigma(\rho)\lbrack V\setminus Y\rbrack)\geq\sum_{i,j}b_{ij}$
(the proof goes along the lines of \cite[Sec.\ 7--9]{Knuth}).
Hence,
	$$\thetI(\sigma(\rho))\geq\thetI(\sigma(\rho)\lbrack V\setminus Y\rbrack)\geq
		\sum_{i,j}b_{ij}\geq nr^{-1/2}/(2c_2).$$
Finally, applying L.\ \ref{Lemma_ConcReg} once more yields our assertion.
\qed\end{proof}

It is shown in \cite{Wormald} that in the case $r=o(n^{1/4})$ we have
	$\pr(\pi(\rho)\textrm{ is a simple graph})\geq\exp(-o(n^{1/2})).$
Therefore, letting $t=\Omega(nr^{-1/2})$ and applying L.\ \ref{Lemma_ConcReg} we conclude that
there is some constant $c>0$ such that
	$\pr(\thetII(\gnr)\leq c nr^{-1/2})=1-o(1)$.
A similar argument proves the lower bound on $\thetI(\gnr)$.

\section{Approximating the Independence Number and Deciding $k$-colorability}\label{Sec_Algo}

In this section we present the algorithms required for Thms.\ 6--8.
The algorithm for the independent set problem is essentially identical
with that proposed in \cite{CojaTaraz},
and the algorithm for deciding $k$-colorability resembles that given in \cite{KrivDecide}.
Thus, our contribution is that using our new results on the \Lovasz\ number of random graphs
and the vector chromatic number,
we can improve on the \emph{analyses} given in \cite{CojaTaraz,KrivDecide}.

\subsubsection{Approximating the independence number.}
The algorithm \texttt{ApproxMIS} for approximating the independence number
consists of two parts.
First, we employ a certain greedy procedure that on input $G=\gnp$ most
probably finds an independent set of size at least $\ln(np)/(2p)$, thereby
providing a lower bound on $\alpha(G)$.
Secondly, we compute $\thet(G)$ to bound $\alpha(G)$ from above.
Throughout, we assume that $np\geq c_0$ for some large constant $c_0$.
%If $\thet(G)$ is not much larger than its probable value (cf.\ Thm.\ \ref{Thm_thetgnp}),
%then the independent set computed by the
%greedy procedure is within the desired approximation ratio and we halt.
%Otherwise, \texttt{ApproxMIS} will spend superpolynomial time on input $G$
%to ensure the approximation ratio stated in Thm.\ \ref{Thm_MIS}.

Following \cite{KrivVu}, to find a large independent set of $G=\gnp$, we run the greedy algorithm
for graph coloring and pick the largest color class it produces.
Remember that the greedy algorithm goes through the vertices $v=1,\ldots,n$
of $G$, and assigns to $v$ the least color among $\{1,\ldots,n\}$
that is not occupied by a neighbor $w<v$ of $v$.

\begin{lemma}\label{Lemma_GreedyMIS}
The probability that the largest color class produced by the greedy coloring
algorithm contains $<\ln(np)/(2p)$ vertices is at most $\exp(-n)$.
\end{lemma}
\begin{proof}
The proof given in \cite{KrivVu} for the case that $p\geq n^{\eps-1/2}$ carries over.
\qed\end{proof}

The following algorithm is essentially identical with the one given in \cite{CojaTaraz}.
(The difference between the algorithm proposed in \cite{KrivVu} and the one below
is that our algorithm uses the \Lovasz\  number as an upper bound on $\alpha(G)$
instead of the largest eigenvalue of the matrix (\ref{eqdefMstrich}).)

\begin{Algo}\label{Alg_ApproxMIS}\upshape\texttt{ApproxMIS$(G)$}\\
\emph{Input:} A graph $G=(V,E)$.\
\emph{Output:} An independent set of $G$.
\begin{enumerate}
\item Run the greedy algorithm for graph coloring on input $G$.
  Let $I$ be the largest resulting color class.
  If $\#I<\ln(np)/(2p)$, then go to 5.
\item Compute $\thet(G)$.
	If $\thet(G)\leq C\sqrt{n/p}$, then output $I$ and terminate.
	Here $C$ denotes some sufficiently large constant (cf.\ the analysis below).
\item Check whether there exists a subset $S$ of $V$, $\#S=25\ln(np)/p$, such
	that $\#V\setminus(S\cup N(S))>12(n/p)^{1/2}$.
	If no such set exists, then output $I$ and terminate.
\item Check whether in $G$ there is an independent set of size $12(n/p)^{1/2}$.
	If this is not the case, then output $I$ and terminate.
\item Enumerate all subsets of $V$ and output a maximum independent set.
\end{enumerate}
\end{Algo}

\begin{lemma}\label{Lemma_timeMIS}
The expected running time of \texttt{ApproxMIS}$(\gnp)$ is polynomial.
\end{lemma}
\begin{proof}
The first two steps can be implemented in polynomial time,
because $\thet(G)$ can be computed efficiently \cite{GLS} (we may disregard
rounding issues in this paper).
By Thm.\ \ref{Thm_thetgnp}, the median $\mu$ of $\thet(\gnp)$ is at most
$c\sqrt{n/p}$, for some constant $c$.
Therefore, Thm.\ \ref{Thm_Largedev} entails that the probability that
\texttt{ApproxMIS} runs step 3 is less than $\exp(-(n/p)^{1/2})$, provided
$C$ is large enough.
Furthermore, up to polynomial factors, step 3 consumes time
	$\leq\exp(25\ln(np)^2/p)<\exp(\sqrt{n/p}).$
Hence, the expected time spent executing step 3 is polynomial.
Taking into account L.\ \ref{Lemma_GreedyMIS},
the expected running time of the remaining steps can be estimated
as in the proof of Thm.\ 4 in \cite{CojaTaraz}.
\qed\end{proof}

Finally, we claim that there exists some constant $\zeta>0$ such that
\texttt{ApproxMIS}$(G)$ finds an independent set of
size at least $\alpha(G)\ln(np)/(\zeta\sqrt{np})$, for all graphs $G$.
Since step 5 will always find an independent set of size $\alpha(G)$,
we may assume that $\#I\geq\ln(np)/(2p)$.
Thus, if $\alpha(G)\leq C'\sqrt{n/p}$, then
	$\#I/\alpha(G)\geq\ln(np)/(2C'\sqrt{np}).$
Further, it is easily seen that in the case $\alpha(G)> C'\sqrt{n/p}$ the algorithm
\texttt{ApproxMIS} will run step 5.

\begin{remark}
The lower bounds on $\thetI(\gnp)$ in Thm.\ \ref{Thm_thetgnp} shows that
we could not achieve an approximation ratio of $(np)^{\beta}$, $\beta<1/2$,
even if we would use the relaxation $\thetI$ instead of $\thet$
to upper-bound the independence number in our algorithm \texttt{ApproxMIS}.
\end{remark}

\subsubsection{Deciding $k$-colorability.}

Following \cite{KrivDecide},
we decide $k$-colorability by computing the vector chromatic number of the input graph.
Let $k=k(n)$ be a sequence of positive integers such that $k(n)=o(\sqrt{n})$.
Since the vector chromatic number is always a lower bound on the chromatic number,
the answer of the following algorithm is correct for all input graphs $G$.

\begin{Algo}\label{Alg_Decide}\upshape$\texttt{Decide}_k(G)$\\
\emph{Input:} A graph $G=(V,E)$.\
\emph{Output:} Either ``$\chi(G)\leq k$'' or ``$\chi(G)>k$''.
\begin{enumerate}
\item If $\vchi(G)>k$ then terminate with output ``$\chi(G)>k$''.
\item Otherwise, compute $\chi(G)$ in time $o(\exp(n))$ using
	Lawler's algorithm \cite{Lawler}, and answer correctly.
\end{enumerate}
\end{Algo}

Thm.\ \ref{Thm_Decide} is a consequence of the following lemma.

\begin{lemma}\label{Lemma_timeDecide}
Suppose that $p\geq Ck^2/n$ for some large constant $C$.
Then the expected running time of $\texttt{Decide}_k(\gnp^+)$ is polynomial.
\end{lemma}
\begin{proof}
In \cite{KMS} it is shown that $\bthetI$ can be computed in polynomial time
(we disregard rounding issues).
Since the second step consumes time $o(\exp(n))$, inequality (\ref{eqexp})
shows that the expected running time of $\texttt{Decide}_k$ on input $\gnp$
is polynomial.
Consequently, by (\ref{eqmonotone}),
we conclude that the expected running time of $\texttt{Decide}_k(\gnp^+)$ is polynomial.
\qed\end{proof}

The analysis of $\texttt{Decide}_k$ on input $\gnr$, $r\geq Ck^2$,
is based on (\ref{eqvchignr}) and yields the proof of Thm.\ \ref{Thm_DecideReg}.

\begin{remark}
The analysis of $\texttt{Decide}_k$ shows that we can decide in polynomial expected time
whether $\gnp^+$ is $k$-colorable, provided $np\geq c_0k^2$.
Conversely, the upper bounds on $\bthetI(\gnp)$, $\bthet(\gnp)$, $\bthetII(\gnp)$
in Thm.\ \ref{Thm_thetgnp} show that even if we would use the relaxation $\bthetII$
instead of the vector chromatic number $\bthetI$ in our algorithm
$\texttt{Decide}_k$, we would still have to assume that $np=\Omega(k^2)$.
\end{remark}

\section{Conclusion}

The results presented in this paper
show the \Lovasz\ number and other SDP relaxations of the independence number
or the chromatic number provide powerful tools in the design of algorithms
with a polynomial expected running time.
Indeed, in addition to the algorithmic applications given in this paper,
the results were used by Coja-Oghlan, Goerdt, Lanka, and Sch\"adlich
to obtain an algorithm for deciding in polynomial expected time whether a random $2k$-SAT formula is
satisfiable \cite{GoerdtCoja}.
Furthermore, using Thms.\ \ref{Thm_Largedev} and \ref{Thm_thetgnp},
Coja-Oghlan obtained an algorithm that finds
a large independent set hidden in a semirandom graph in polynomial expected time~\cite{FindLarge}.
Therefore, the author expects that the general idea of combining large deviation techniques
such as Talagrand's inequality with SDP relaxations will lead to further
contributions to the algorithmic theory of random structures (cf.\ also \cite{MAXCUT}).

In comparison with purely combinatorial techniques or computing eigenvalues,
semidefinite programming requires a rather heavy machinery (cf.\ \cite{GLS}).
However, compared to the eigenvalues of the adjacency matrix, semidefinite
programs such as the \Lovasz\ number seem to be rather ``robust'' (cf.\ the
discussion in \cite{Feige}).
In fact, this robustness may be the reason why for the \Lovasz\ number we can derive
a large deviation result such as Thm.\ \ref{Thm_Largedev}.
By contrast, for ``small'' values of $p$
no similarly strong tail bounds on the eigenvalues of the auxiliary
matrix $M(\gnp)$ are known (cf.~\cite{AlonKrivVu,KrivVu})--although
we used to the eigenvalues of $M(\gnp)$ to bound the mean of $\thet(\gnp)$.

\subsubsection{Acknowledgment.}
I am grateful to M.\ Krivelevich and C.\ Helmberg for helpful discussions, and to
U.\ Feige and E.\ Ofek for providing me with their technical report \cite{Ofek}.

\end{document}